# A wavelet integral collocation method for nonlinear boundary value problems

Lei Zhang, Jizeng Wang[1], Xiaojing Liu, Youhe Zhou

*Key Laboratory of Mechanics on Disaster and Environment in Western China, the Ministry of Education, College of Civil Engineering and Mechanics, Lanzhou University, Lanzhou 730000,China*

**Abstract:** A high-order wavelet integral collocation method (WICM) is developed for general nonlinear boundary value problems. This method is established based on Coiflet approximation of multiple integrals of interval bounded functions combined with an accurate and adjustable boundary extension technique. The convergence order of this approximation has been proven to be $N$ as long as the Coiflet with $N$-1 vanishing moment is adopted, which can be any positive even integers. Before the conventional collocation method is applied to the general problems, the original differential equation is changed into its equivalent form by denoting derivatives of the unknown function as new functions and constructing relations between the low- and high-order derivatives. For the linear cases, error analysis has proven that the proposed WICM is order $N$, and condition numbers of relevant matrices are almost independent of the number of collocation points. Numerical examples of a wide range of nonlinear differential equations demonstrate that accuracy of the proposed WICM is even greater than $N$, and most interestingly, such accuracy is independent of the order of the differential equation to be solved. Comparison to existing numerical methods further justifies the accuracy and efficiency of the proposed method.

**Key Words:** nonlinear boundary value problems; Coiflet; collocation method; high-order accuracy

## 1. Introduction

Nonlinear boundary value problems (BVPs) [1] arise from almost every scientific and engineering field, especially mechanical theory of beam and plate structures [2-4]. Commonly, the nonlinear BVPs can be written into a general form

---

[1] Corresponding author, Email: jzwang@lzu.edu.cn

[5]

$$\mathbf{T}\left(x, y, \frac{dy}{dx}, \cdots, \frac{d^n y}{dx^n}\right) = 0, \quad x \in [a,b] \tag{1}$$

with corresponding boundary conditions, where $\mathbf{T}$ is a nonlinear operator and $y$ is an unknown function of $x$.

Since the solutions of nonlinear BVPs are critically important for analyzing scientific and engineering problems, finding accurate and efficient methods for solving Eq. (1) has attracted considerable research attention. In the last few decades, numerous methods have been developed [6-12]: for example, Wang and Wu [7] proposed a fourth-order compact finite difference method (CFDM) to solve nonlinear $2n$th-order multi-point boundary value problems; Geng [8] studied the nonlinear four-point boundary value problems by applying the reproducing kernel Hilbert space methods (RKHSM); Behroozifar [10] suggested a spectral method based on Bernstein polynomials (SMBP) for nonlinear differential equations with multi-point boundary conditions [10]; Liu et al. [6, 12] proposed a wavelet Galerkin method (WGM) for studying nonlinear differential equations with Dirichlet and Neumann boundary conditions. Moreover, there are still many other solution methods including the shooting method, series method, function space method [13], homotopy analysis method [14] and high order finite difference method [15] etc. Although, there are many alternative ways to study the nonlinear BVPs [1], yet finding highly accurate solutions to general nonlinear BVPs remains a challenge.

Wavelet theory is a newly developed powerful mathematical tool that has shown its potential in the numerical analysis of differential equations [16, 17]. Wavelet-based methods combined with conventional Galerkin [16, 18] or collocation techniques [17, 19] have been proposed for various nonlinear engineering problems. In our recent works [6, 20-22], an efficient WGM has been proposed to solve the Bratu equation [6], large deformation bending of circular and rectangular plates [20, 22], and large deflection and post-buckling analysis of nonlinearly elastic rods [21]. The numerical results have considerably better accuracy than many other numerical methods, and show their applicability to strong nonlinear problems. However, just like most other numerical methods for the nonlinear BVPs, the convergence rate of the WGM can be affected by the highest derivatives involved in the equations [4, 6, 20, 21].

In this study, we construct a modified wavelet approximation of $n$-tuple integrals of a function, whose accuracy is found higher than the approximations of its derivatives, and interestingly independent of the tuple of the integral. Based on this finding, we transform the original nonlinear BVPs into its equivalent forms by defining the derivatives of unknown functions as new functions and constructing relations between the low- and high- order derivatives. When solving the transformed equation by using conventional collocation method in terms of the proposed wavelet approximation on $n$-tuple integrals, we show that the accuracy order can be any positive even integer $N=4, 6, 8, \ldots$, for the Coiflet with compact support [3$N$-1]. And this accuracy order is independent of the orders of the original differential equations. Both error analysis and nontrivial numerical examples are given for justifications.

## 2. Wavelet approximation of multiple integrals

The multi-resolution analysis of wavelet theory [23] states that the function space $\mathbf{L}^2(\mathbf{R})$ can be divided into a sequence of nested subspaces $\{0\} \cdots \subset V_0 \subset V_1 \subset \cdots \subset V_j \subset V_{j+1} \subset \cdots \subset \mathbf{L}^2(\mathbf{R})$. A set of orthogonal basis of the subspace $V_j$ can be formed by

$$\phi_{j,k}(x) = 2^{\frac{j}{2}} \phi(2^j x - k), \quad k \in \mathbf{Z}, \tag{2}$$

where $\phi(x)$ is the so-called orthogonal scaling function. A function $f(x) \in \mathbf{L}^2(\mathbf{R})$ can be approximated by projecting this function from $\mathbf{L}^2(\mathbf{R})$ to $V_j$ as

$$f(x) \approx \mathbf{P}^j f(x) = \sum_k c_{j,k} \phi_{j,k}(x) \tag{3}$$

where $c_{j,k} = \int_{-\infty}^{\infty} f(x) \phi_{j,k}(x) dx$, and integer $j$ is the so-called resolution level. The scaling function with compact support can be constructed by using a finite number of low-pass filter coefficients $p_k$ in terms of the relation below:

$$\phi(x) = \sum_k p_k \phi(2x - k) \tag{4}$$

in which subscript $k = 0, 1, 2, \ldots, 3N$-1 for the Coiflet-type wavelet, and $N$-1 is the number of vanishing moments of the corresponding wavelet function [24]. Such a scaling function has the unique property of shifted vanishing moments:

$$\int_{-\infty}^{\infty} (t - M_1)^k \phi(t) dt = 0, \quad 1 \leq k < N \tag{5}$$

where $M_1 \triangleq \int_{-\infty}^{\infty} x \phi(x) dx = \sum_{k \in \mathbf{Z}} p_k k / 2$ is the first-order moment of the scaling function. Based on this unique property, one has $c_{j,k} \approx 2^{-j/2} f(\frac{k+M_1}{2^j})$, such that the approximation of the function can be written as [25]

$$f(x) \approx \mathbf{P}^j f(x) \approx \tilde{\mathbf{P}}^j f(x) = \sum_{k=-\infty}^{\infty} f_{M_1+k} \phi(2^j x - k) \tag{6}$$

where $f_{M_1+k} = f(x_{M_1+k})$, $x_{M_1+k} = (M_1+k)/2^j$, and $f(x) = \mathbf{P}^j f(x) = \tilde{\mathbf{P}}^j f(x)$ when $f(x)$ is any polynomial with an order up to $N$-1, and $\left\| f(x) - \tilde{\mathbf{P}}^j f(x) \right\|_{\mathbf{L}^2(\mathbf{R})} = O(2^{-jN})$ as long as $f(x) \in \mathbf{L}^2(\mathbf{R}) \cap \mathbf{C}^N(\mathbf{R})$ [26].

If the function $f(x)$ is defined on an interval, for instance, [0, 1], then Eq. (6) can be rewritten as

$$f(x) \approx \tilde{\mathbf{P}}^j f(x) = \sum_{k=2-3N}^{2^j-1} f_{k+M_1} \phi(2^j x - k). \tag{7}$$

We define the *n*-tuple integral of the function $f(x)$ as [27, 28]

$$f^{\int_n}(x) \triangleq \int_0^x \int_0^{\xi_n} \cdots \int_0^{\xi_2} f(\xi_1) d\xi_1 d\xi_2 \cdots d\xi_n. \tag{8}$$

Substituting Eq. (7) to Eq. (8) yields

$$f^{\int_n}(x) \approx f_{\tilde{\mathbf{P}}^j}^{\int_n}(x) = \sum_{k=2-3N}^{2^j-1} f_{k+M_1} \phi_{j,k}^{\int_n}(x), \tag{9}$$

where

$$\phi_{j,k}^{\int_n}(x) \triangleq \int_0^x \int_0^{\xi_n} \cdots \int_0^{\xi_2} \phi(2^j \xi_1 - k) d\xi_1 d\xi_2 \cdots d\xi_n,$$

$$f_{\tilde{\mathbf{P}}^j}^{\int_n}(x) = \int_0^x \int_0^{\xi_n} \cdots \int_0^{\xi_2} f_{\tilde{\mathbf{P}}^j}^{\int_n}(\xi_1) d\xi_1 d\xi_2 \cdots d\xi_n.$$

**Theorem 1.** If $f(x) \in \mathbf{L}^2(\mathbf{R}) \cap \mathbf{C}^N(\mathbf{R})$, then the accuracy of approximation (9) can be estimated as

$$\left\| f^{\int_n}(x) - f_{\tilde{\mathbf{P}}^j}^{\int_n}(x) \right\|_{\mathbf{L}^2[0,1]} \leq 2^{-jN} \frac{\Omega(3N-1)}{N!(n-1)!} \left| (M_1-\tau)^N \Theta \phi(\tau) \right| \cdot \left\| x^{n-1} \right\|_{\mathbf{L}^2[0,1]} \tag{10}$$

and

$$| f^{\int_n} - f_{\tilde{\mathbf{P}}^j}^{\int_n} | \leq 2^{-jN} \frac{\Omega(3N-1)}{N!(n-1)!} |(M_1-\tau)^N \Theta \varphi(\tau) x^{n-1}| \tag{11}$$

where $\Theta = \max[| f^{(N)}(x) |], x \in [0, 3N-1]$, $\Omega = 1 + (3N-2)/2^j$, and $\tau \in [0, 3N-1]$.

**Proof.** Using Taylor expansion of $f(y)$ at the point $x$ gives:

$$f(y) = \sum_{n=0}^{N-1} \left[ \frac{f^{(n)}(x)}{n!} (y-x)^n \right] + \frac{f^{(N)}(\vartheta)}{N!} (y-x)^N \qquad (12)$$

where $f^{(n)}(x) \triangleq d^n f(x)/dx^n$ and $\vartheta$ is on the segment connecting $y$ and $x$. Assigning $y = (k+M_1)/2^j$, $k \in \mathbf{Z}$ into Eq. (12) gives the expansion of $f_{k+M_1}$ as

$$f_{k+M_1} = \sum_{n=0}^{N-1} \left[ \frac{f^{(n)}(x)}{n!} (\frac{k+M_1}{2^j} - x)^n \right] + \frac{f^{(N)}(\vartheta)}{N!} (\frac{k+M_1}{2^j} - x)^N .$$

Further considering the property of Coiflet expansion $\sum_{k \in \mathbf{Z}} (k+M_1 - 2^j x)^n \phi(2^j x - k) = 0^n$ for $0 \leq n < N$, then for $x \in [0,1]$, Eq. (7) can be rewritten into [26]

$$\tilde{\mathbf{P}}^j f(x) = f(x) + 2^{-jN} \sum_{2^j > k > 2-3N} \frac{f^{(N)}(\vartheta_{j,k})}{N!} (k+M_1 - 2^j x)^N \phi(2^j x - k), \qquad (13)$$

where $\vartheta_{j,k}$ becomes locating between $y = (k+M_1)/2^j$ and $x$, and $\mathrm{Supp}[\phi(x)] = [0, 3N-1]$ is considered for determining the range of summation index $k$. The $n$-tuple integral of Eq. (13) can be expressed as

$$f_{\tilde{\mathbf{P}}^j}^{\int_n}(x) = f^{\int_n}(x) + 2^{-jN} \sum_{2^j > k > 2-3N} \int_0^x \int_0^{\xi_n} \cdots \int_0^{\xi_2} \frac{f^{(N)}(\vartheta_{j,k})}{N!} (k+M_1 - 2^j \xi_1)^N \phi(2^j \xi_1 - k) d\xi_1 d\xi_2 \cdots d\xi_n . \qquad (14)$$

In addition, denoting $\xi_1' = 2^j \xi_1 - k$, we have

$$\begin{aligned} f_{\tilde{\mathbf{P}}^j}^{\int_n}(x) &= f^{\int_n}(x) + 2^{-j(N+1)} \sum_{2^j > k > 2-3N} \int_0^x \int_0^{\xi_n} \cdots \int_{-k}^{2^j \xi_2 - k} \frac{f^{(N)}(\vartheta_{j,k})}{N!} (M_1 - \xi_1')^N \phi(\xi_1') d\xi_1' d\xi_2 \cdots d\xi_n \\ &= f^{\int_n}(x) + 2^{-j(N+1)} \int_0^x \int_0^{\xi_n} \cdots \int_0^{\xi_3} g(\xi_2) d\xi_2 \cdots d\xi_n \end{aligned} \qquad (15)$$

where $g(\xi_2)$ is defined as

$$g(\xi_2) = \sum_{2^j > k > 2-3N} \int_{-k}^{2^j \xi_2 - k} \frac{f^{(N)}(\vartheta_{j,k})}{N!} (M_1 - \xi_1')^N \phi(\xi_1') d\xi_1' . \qquad (16)$$

Considering again $\mathrm{Supp}[\phi(x)] = [0, 3N-1]$, then $g(\xi_2)$ can be further written into

$$g(\xi_2) = \sum_{2^j > k > 2-3N} \int_{[0,3N-1] \cap [-k, 2^j \xi_2 - k]} \frac{f^{(N)}(\vartheta_{j,k})}{N!} (M_1 - \xi_1')^N \phi(\xi_1') d\xi_1' , \qquad (17)$$

so that

$$\begin{aligned} |g(\xi_2)| &= \left| \sum_{2^j > k > 2-3N} \int_{[0,3N-1] \cap [-k, 2^j \xi_2 - k]} \frac{f^{(N)}(\vartheta_{j,k})}{N!} (M_1 - \xi_1')^N \phi(\xi_1') d\xi_1' \right| \\ &\leq (2^j + 3N - 2)(3N - 1) \left| \frac{\Theta}{N!} (M_1 - \tau)^N \phi(\tau) \right| \end{aligned} \qquad (18)$$

where $\Theta = \max[|f^{(N)}(x)|]$, $x \in [0, 3N-1]$, and $\tau \in [0, 3N-1]$. Then the accuracy of Eq. (14) can be estimated as

$$\left\| f^{\int_n}(x) - f^{\int_n}_{\tilde{\mathbf{P}}^j}(x) \right\|_{\mathbf{L}^2[0,1]} = 2^{-j(N+1)} \left\| \int_0^x \int_0^{\xi_n} \cdots \int_0^{\xi_3} g(\xi_2) d\xi_2 \cdots d\xi_n \right\|_{\mathbf{L}^2[0,1]}$$

$$\leq 2^{-j(N+1)} \frac{(2^j + 3N - 2)(3N-1)|\Theta(M_1-\tau)^N \phi(\tau)|}{N!} \left\| \int_0^x \int_0^{\xi_n} \cdots \int_0^{\xi_3} d\xi_2 \cdots d\xi_n \right\|_{\mathbf{L}^2[0,1]} \quad (19)$$

$$\leq 2^{-j(N+1)} \frac{(2^j + 3N - 2)(3N-1)}{N!(n-1)!} |(M_1-\tau)^N \Theta \phi(\tau)| \cdot \|x^{n-1}\|_{\mathbf{L}^2[0,1]}$$

$$\leq 2^{-jN} \frac{\Omega(3N-1)}{N!(n-1)!} |(M_1-\tau)^N \Theta \phi(\tau)| \cdot \|x^{n-1}\|_{\mathbf{L}^2[0,1]}$$

where $\Omega = 1 + (3N-2)/2^j$. Similar to Eq. (19), we also have

$$| f^{\int_n} - f^{\int_n}_{\tilde{\mathbf{P}}^j} | \leq 2^{-jN} \frac{\Omega(3N-1)}{N!(n-1)!} |(M_1-\tau)^N \Theta \varphi(\tau) x^{n-1}| \quad (20)$$

In the case that $j \geq 3$, and Coiflet with $N = 6$ and $M_1 = 7$ is adopted, whose filter coefficients $p_k$ are listed in Table 1, we have $1 \leq \Omega \leq 3$.

The procedure for calculating $\phi^{\int_n}_{j,k-M_1}(x)$ defined by Eq. (9) can be found in Appendix A. In all the following numerical examples, Coiflet with $N = 6$ and $M_1 = 7$ will be adopted, corresponding values of $\phi^{\int_n}(k)$ for $k = 1, 2, \ldots, 17$ and $n = 1, 2, 3, 4$ have been given in Table 2. To numerically demonstrate the accuracy of the approximation (9), Fig. 1 shows the absolute error of $n$-tuple integrals, $\int_0^1 \int_0^{\xi_n} \cdots \int_0^{\xi_2} \sin(\pi \xi_1) d\xi_1 d\xi_2 \cdots d\xi_n$, for $n = 1, 2, 3$, and 4. It can be seen from Fig. 1 that the absolute error of these integrals decreases with the resolution level $j$, indicating a convergence rate of order 7.5 which is even higher than $N=6$ that shows in Eq. (11).

When considering an interval bounded function, for example, $f(x)$ with $x \in [0, 1]$, we note that $x_k = (M_1+k)/2^j$ in Eq. (9) can locate outside the interval $[0, 1]$, where $f(x)$ may have no definitions. Therefore, the techniques of boundary extension are usually adopted to resolve this problem. In this study, the power series expansion is used at each boundary as [6]

$$f(x) = \begin{cases} \sum_{i=0}^{N-1} \dfrac{F_{0,i} x^i}{i!} & x < 0 \\ f(x) & 0 \leq x \leq 1 \\ \sum_{i=0}^{N-1} \dfrac{F_{1,i}(x-1)^i}{i!} & x > 1 \end{cases}, \quad (21)$$

where coefficients $F_{0,i}$ and $F_{1,i}$ can be given by

$$F_{1,i} = 2^{ij} \sum_{k=0}^{\alpha_2} \zeta_{1,i,k} f_{2^j-k}, \tag{22}$$

and

$$F_{0,i} = 2^{ij} \sum_{k=0}^{\alpha_1} \zeta_{0,i,k} f_k. \tag{23}$$

To derive the coefficients $\zeta_{1,i,k}$, substituting Eq. (22) into Eq. (21), then into Eq. (7), taking the derivatives on both sides of the resulting equation with respect to $x$, and considering $x = 1$, we have

$$\frac{d^i f(1)}{dx^i} \approx 2^{ij} \sum_{k=2^j-3N+2+M_1}^{2^j} f_k \phi^{(i)}(2^j - k + M_1) + 2^{ij} \sum_{k=2^j+1}^{2^j-1+M_1} \sum_{l=0}^{N-1} \frac{F_{1,l}}{i!} \left(\frac{k}{2^j} - 1\right)^l \phi^{(i)}(2^j - k + M_1). \tag{23}$$

where $j$ is assumed to be sufficiently large so that $2^j - 3N + 2 + M_1 > 0$. By using $F_{1,i}$ to replace $d^i f(1)/dx^i$, inserting $F_{1,i} = 2^{ij} \sum_{k=0}^{\alpha_2} \zeta_{1,i,k} f_{2^j-k}$ into Eq. (24), and taking $\alpha_2 = 3N - 2 - M_1$, we have

$$\sum_{k=0}^{\alpha_2} [\phi^{(i)}(k + M_1) - \zeta_{1,i,k} + \sum_{l=0}^{N-1} \zeta_{1,l,k} \sum_{t=1}^{M_1-1} \frac{1}{l!} t^l \phi^{(i)}(M_1 - t)] f_{2^j-k} = 0. \tag{24}$$

Eq. (25) can be written in a matrix form as

$$[\mathbf{R}_1 - (\mathbf{I} - \mathbf{S}_1)\mathbf{P}_1]\mathbf{F} = 0, \tag{25}$$

where $\mathbf{F} = \{f_{2^j-k}\}$, $\mathbf{P}_1 = \{\zeta_{1,i,k}\}$, $\mathbf{R}_1 = \{\phi^{(i)}(k + M_1)\}$, $\mathbf{S}_1 = \{\sum_{t=1}^{M_1-1} \frac{1}{l!} t^l \phi^{(i)}(M_1 - t)\}$, $i, l = 0, 1, 2, ... N - 1$, and $k = 0, 1, 2, ... \alpha_2$. Forcing Eq. (26) to be satisfied for any $\mathbf{F}$, then from Eq. (26), we obtain

$$\mathbf{P}_1 = (\mathbf{I} - \mathbf{S}_1)^{-1} \mathbf{R}_1. \tag{26}$$

Similarly, we can obtain

$$\mathbf{P}_0 = (\mathbf{I} - \mathbf{S}_0)^{-1} \mathbf{R}_0, \tag{27}$$

where the matrices are $\mathbf{P}_0 = \{\zeta_{0,i,k}\}$, $\mathbf{R}_0 = \{\phi^{(i)}(M_1 - k)\}$, $\mathbf{S}_0 = \{\sum_{t=-\alpha_1}^{-1} \frac{1}{l!} t^l \phi^{(i)}(M_1 - t)\}$, $\alpha_1 = M_1 - 1$, $i, l = 0, 1, 2, ... N - 1$, and $k = 0, 1, 2, ... \alpha_1$. Inserting Eq. (21) into Eq. (9) and making further rearrangements as shown in [6] yields

$$f^{j_n}(x) \approx \sum_{k=0}^{2^j} f_k \Phi_{j,k}^{j_n}(x), \tag{28}$$

where

$$\Phi_{j,k}^{\int_n}(x) = \begin{cases} \phi_{j,k-M_1}^{\int_n}(x) + \sum_{m=-\alpha_2}^{0} T_{L,k}(\frac{m}{2^j})\phi_{j,m-M_1}^{\int_n}(x), & 0 \le k \le \alpha_1 \\ \phi_{j,k-M_1}^{\int_n}(x), & \alpha_1+1 \le k \le 2^j-\alpha_2-1 \\ \phi_{j,k-M_1}^{\int_n}(x) + \sum_{m=1}^{\alpha_1} T_{R,2^j-k}(1+\frac{m}{2^j})\phi_{j,2^j+m-M_1}^{\int_n}(x), & 2^j-\alpha_2 \le k \le 2^j \end{cases}, \quad (29)$$

and

$$T_{L,l}(x) = \sum_{i=0}^{N-1} 2^{ij} \frac{\zeta_{a,i,l}}{i!} x^i, \quad T_{R,l}(x) = \sum_{i=0}^{N-1} 2^{ij} \frac{\zeta_{b,i,l}}{i!} (x-1)^i. \quad (30)$$

Comparing to Ref. [6], Eq. (29) and (30) provide an accurate boundary extension technique. We note that the boundary extension guaranteed the ($N$-1)th order smoothness of a function near the domain boundaries [26]. Although such an operation is expected to effectively reduce $\Theta$ in Eq. (10), yet the order of the accuracy of Eq. (29) should still be estimated by Eq. (10).

## 3. Wavelet integral collocation method (WICM)

In order to establish a wavelet-based numerical method for the solution of the general nonlinear BVPs as shown in (1), we assume that $x \in [0,1]$ and the boundary conditions located at the end points $x = 0, 1$. Defining various derivatives of the unknown function, $d^i y / dx^i$, as new functions $y_i$, converts Eq. (1) into

$$\begin{cases} \mathbf{T}(x, y_0, y_1, \cdots, y_n) = 0 \\ y_0 = \int_0^x \int_0^{\xi_n} \cdots \int_0^{\xi_2} y_n(\xi_1) d\xi_1 d\xi_2 \cdots d\xi_n + \sum_{i=0}^{n-1} \frac{x^i}{i!} y_i(0) \\ y_1 = \int_0^x \int_0^{\xi_{n-1}} \cdots \int_0^{\xi_2} y_n(\xi_1) d\xi_1 d\xi_2 \cdots d\xi_{n-1} + \sum_{i=1}^{n-1} \frac{x^{i-1}}{(i-1)!} y_i(0), \quad x \in [0,1], \quad (31) \\ \vdots \\ y_{n-1} = \int_0^x y_n d\xi + y_{n-1}(0) \end{cases}$$

where $y_i(0)$ can be determined by the boundary conditions at the end points.

We then approximate the newly defined function, $y_n$, by

$$y_n \approx \sum_{k=0}^{2^j} y_n(\frac{k}{2^j}) \Phi_{j,k}(x), \quad (32)$$

Substituting Eq. (33) into Eq. (32) leads to

$$\begin{cases} \mathbf{T}(x, y_0, y_1, \cdots, y_n) = 0 \\ y_0(x) \approx \sum_{k=0}^{2^j} y_{n,k} \Phi_{j,k}^{\int_n}(x) + \sum_{i=0}^{n-1} \frac{x^i}{i!} y_i(0) \\ y_1(x) \approx \sum_{k=0}^{2^j} y_{n,k} \Phi_{j,k}^{\int_{n-1}}(x) + \sum_{i=1}^{n-1} \frac{x^{i-1}}{(i-1)!} y_i(0) \\ \vdots \\ y_{n-1}(x) \approx \sum_{k=0}^{2^j} y_{n,k} \Phi_{j,k}^{\int}(x) + y_{n-1}(0) \end{cases} \quad (33)$$

The conventional collocation method is applied to Eq. (34), then we obtain

$$\begin{cases} \mathbf{T}(x_l, y_{0,l}, y_{1,l}, \cdots, y_{n,l}) = 0 \\ y_{0,l} \approx \sum_{k=0}^{2^j} y_{n,k} \Phi_{j,k}^{\int_n}(x_l) + \sum_{i=0}^{n-1} \frac{x_l^i}{i!} y_{i,0} \\ y_{1,l} \approx \sum_{k=0}^{2^j} y_{n,k} \Phi_{j,k}^{\int_{n-1}}(x_l) + \sum_{i=1}^{n-1} \frac{x_l^{i-1}}{(i-1)!} y_{i,0} \\ \vdots \\ y_{n-1,l} \approx \sum_{k=0}^{2^j} y_{n,k} \Phi_{j,k}^{\int}(x_l) + y_{n-1,0} \end{cases} \quad (34)$$

where $x_l = l/2^j$ and $y_{i,l} = y_i(x_l)$, $l = 0,1,2,\cdots,2^j$. The boundary values of $y_i(0)$ that are not explicitly given usually can be determined from boundary conditions at other points, which can be eventually expressed in terms of the values of unknown function $y_n$ and the integral of scaling function at nodal points $x_l$. This procedure can be explicitly demonstrated in the next subsection by applying this method to a general fourth order linear differential equation.

Eq. (35) can be rewritten into a matrix form as

$$\mathbf{T}(\mathbf{x}, \mathbf{A}_n \mathbf{y}_n + \mathbf{b}_0, \mathbf{A}_{n-1} \mathbf{y}_n + \mathbf{b}_1, \cdots, \mathbf{A}_1 \mathbf{y}_n + \mathbf{b}_{n-1}) \approx 0, \quad (35)$$

where matrices and vectors are $\mathbf{A}_i = \{a_{i,l,k} = \Phi_{j,k}^{\int_i}(x_l)\}$, $\mathbf{b}_i = \{b_{i,l} = \sum_{k=i}^{n-1} x_l^k y_k(0)/i!\}$, $l$, $k=0, 1, \ldots, 2^j$. We note that the matrix entries $a_{i,l,k}$ are given by the values of $i$-tuple integrals of the scaling function, which can be obtained exactly by using the procedure suggested in Appendix A and the expression of the modified scaling basis given in Eq. (30). The corresponding dimension of the nonlinear algebraic equations in Eq. (36) is $2^j + 1$. Once Eq. (36) is solved, then we can obtain $y_n(x_l)$. The value of $y_i(x_l)$, for $i < n$ can be achieved through the relation $\mathbf{y}_i = \mathbf{A}_{n-i} \mathbf{y}_n + \mathbf{b}_i$.

The WICM can be easily extended to solve high dimensional nonlinear BVPs, since the set of scaling bases for high dimensional space can be directly extended by the tensor products of one-dimensional wavelet bases.

## 4. Error analysis for the solution of linear BVPs

As a typical example of linear BVPs, we consider the general forth order linear differential equation as follows

$$\begin{cases} \sum_{i=0}^{4} \alpha_i \dfrac{d^i u}{dx^i} = f(x) \\ u(0)=a_0, u(1)=a_1, d^2u/dx^2|_{x=0}=b_0, d^2u/dx^2|_{x=1}=b_1 \end{cases}. \quad (36)$$

Defining $u_n = d^n u / dx^n$, Eq. (37) can be transformed into

$$\alpha_4 u_4 + \sum_{i=0}^{3} \alpha_i u_4^{I_{4-i}} + \sum_{i=1}^{4} \alpha_i \sum_{n=0}^{i-1} \dfrac{x^n}{n!} u_n(0) = f(x). \quad (37)$$

To obtain $u_1(0)$ and $u_3(0)$, we consider the following relations

$$u(x) = u_4^{I_4}(x) + \sum_{n=0}^{3} \dfrac{x^n}{n!} u_n(0), \quad (38)$$

$$\dfrac{d^2 u(x)}{dx^2} = u_4^{I_2}(x) + \sum_{n=0}^{1} \dfrac{x^n}{n!} u_n(0). \quad (39)$$

Letting $x=1$, we can obtain

$$a_1 = u_4^{I_4}(1) + a_0 + u_1(0) + \dfrac{b_0}{2} + \dfrac{u_3(0)}{6}, \quad (40)$$

$$b_1 = u_4^{I_2}(1) + a_0 + u_1(0). \quad (41)$$

From Eqs. (41) and (42), we have

$$u_1(0) = b_1 - u_4^{I_2}(1) - a_0, \quad (42)$$

$$u_3(0) = 6a_1 - 6b_1 - 3b_0 - 6u_4^{I_4}(1) + 6u_4^{I_2}(1). \quad (43)$$

Substituting Eqs. (43) and (44) into Eq. (37), and considering the boundary conditions, we have

$$\alpha_4 u_4 + \sum_{i=0}^{3} \alpha_i u_4^{I_{4-i}} + \sum_{i=1}^{4} \alpha_i a_0 + (b_1 - u_4^{I_2}(1) - a_0) \sum_{i=2}^{4} \alpha_i x + b_0 \sum_{i=3}^{4} \alpha_i x^2 / 2 \\ + \alpha_4 [a_1 - b_1 - b_0/2 - u_4^{I_4}(1) + u_4^{I_2}(1)] x^3 = f(x). \quad (44)$$

Eq. (45) can be rewritten into

$$\alpha_4 u_4 + \sum_{i=0}^{3} \alpha_i u_4^{I_{4-i}} - \alpha_4 x^3 u_4^{I_4}(1) - \beta(x) u_4^{I_2}(1) = g(x) \quad (45)$$

where $\beta(x) = \sum_{i=2}^{4} \alpha_i x - \alpha_4 x^3$ and

$$g(x) = f(x) - \sum_{i=1}^{4} \alpha_i a_0 - (b_1 - a_0)\sum_{i=2}^{4} \alpha_i x - \frac{b_0}{2}\sum_{i=3}^{4} \alpha_i x^2 - \alpha_4[a_1 - b_1 - b_0/2]x^3. \quad (46)$$

From Eqs. (11) and (29), we can obtain

$$u_4^{J_{4-i}}(x) = \sum_{k=0}^{2^j} u_{4,k} \Phi_{j,k}^{J_{4-i}}(x) + O(2^{-jN}). \quad (47)$$

Substituting Eq. (48) into Eq. (46) yields

$$\alpha_4 u_4(x) + \sum_{k=0}^{2^j} u_{4,k}\sum_{i=0}^{3}\alpha_i \Phi_{j,k}^{J_{4-i}}(x) - \alpha_4 x^3 \sum_{k=0}^{2^j} u_{4,k}\Phi_{j,k}^{J_4}(1) - \beta(x)\sum_{k=0}^{2^j} u_{4,k}\Phi_{j,k}^{J_2}(1) = g(x) + O(2^{-jN}).(48)$$

Letting $x=0, 1/2^j, ..., 1$, we can have

$$\mathbf{AU} = \mathbf{G} + \mathbf{\Delta} \quad (49)$$

where $\mathbf{U} = \{U_k = u_{4,k}\}^T$, $\mathbf{G} = \{G_l = g(l/2^j)\}^T$, $\mathbf{\Delta} = \{\Delta_l = O(2^{-jN})\}^T$, and

$$\mathbf{A} = \{A_{l,k} = \alpha_4 \delta_{lk} + \sum_{i=0}^{3}\alpha_i \Phi_{j,k}^{J_{4-i}}(\frac{l}{2^j}) - \alpha_4(\frac{l}{2^j})^3 \Phi_{j,k}^{J_4}(1) - \beta(\frac{l}{2^j})\Phi_{j,k}^{J_2}(1)\}, \quad (51)$$

$l, k=0, 1, 2, …, 2^j$. According to Appendix B, we know that $\Phi_{j,k}^{J_i}(\frac{l}{2^j}) = \frac{1}{2^j}K_{i,l,k} = O(\frac{1}{2^j})$,

therefore Eq. (51) can be rewritten as

$$\mathbf{A} = \{A_{l,k} = \alpha_4 \delta_{l,k} + \frac{1}{2^j}M_{l,k}\} = \alpha_4 \mathbf{I} + \frac{1}{2^j}\mathbf{M}$$

where $\mathbf{I}$ is the unit matrix and

$$\mathbf{M} = \{M_{l,k} = \sum_{i=0}^{3}\alpha_i K_{4-i,l,k} - \alpha_4(\frac{l}{2^j})^3 K_{4,2^j,k} - \beta(\frac{l}{2^j})K_{2,2^j,k}]\}, l, k=0, 1, 2, …, 2^j.$$

Obviously, when the resolution level $j$ is large enough, $\mathbf{A} \approx \alpha_4 \mathbf{I}$ becomes a diagonally-dominant matrix, as elements of $\mathbf{M}$ is on the order of $O(1)$ in terms of Appendix B.

Considering Eq. (50), we can obtain the exact solution $\mathbf{U}$ as,

$$\mathbf{U} = \mathbf{A}^{-1}\mathbf{G} + \mathbf{A}^{-1}\mathbf{\Delta} = \bar{\mathbf{U}} + \mathbf{A}^{-1}\mathbf{\Delta} \quad (50)$$

where $\bar{\mathbf{U}} = \mathbf{A}^{-1}\mathbf{G}$ is the approximate solution obtained by using the proposed WICM. From Eq. (52), we have the following error estimation

$$\|\mathbf{U} - \bar{\mathbf{U}}\|_2 = \|\mathbf{A}^{-1}\mathbf{\Delta}\|_2 \leq \|\mathbf{A}^{-1}\|_2 \|\mathbf{\Delta}\|_2 = C2^{-jN}\|\mathbf{A}^{-1}\|_2 \quad (51)$$

in which $C$ is a constant. When the resolution level $j$ is big enough, as we have $\mathbf{A} \approx \alpha_4 \mathbf{I}$, then Eq. (53) can be

$$\|\mathbf{U}-\bar{\mathbf{U}}\|_2 = O(2^{-jN}). \tag{52}$$

Even for small resolution levels, the matrix $\mathbf{A}$ can still have a very good property. For example, we consider three examples with $[\alpha_0, \alpha_1, \alpha_2, \alpha_3, \alpha_4]$ assigned as $[1,1,1,1,1]$, $[1,1,1,1,0.5]$ and $[0.1, 0.1, 0.1, 0.1, 1]$, and matrix $\mathbf{A}$ is correspondingly denoted as $\mathbf{A}_1$, $\mathbf{A}_2$ and $\mathbf{A}_3$. The condition number and the norm of the inverse of these matrices, $K_2(\mathbf{A}_i)$ and $\|\mathbf{A}_i^{-1}\|_2$, $i=1,2,3$, have been given in Table 3. It can be seen that $\|\mathbf{A}_i^{-1}\|_2$ is close to $1/\alpha_4$ and almost independent of the resolution level $j$. The condition number $K_2(\mathbf{A}_i)$ decreases when $j$ increases. In contrary, for most existing numerical methods, the condition number of corresponding stiffness matrix usually significantly increases as the number of grids increases [29, 30], which may be the reason why Bernstein polynomials, although uniformly convergent, are generally considered rather poor for numerical processes [31].

### 4. Numerical examples

In this section, five nonlinear BVPs are studied to demonstrate the convergence and accuracy of the proposed WICM, which includes the typical Bratu equation [6, 10], the nonlinear fourth-order differential equations with various boundary conditions [7–9], the large deflection bending of circular plates [20] governed by coupled nonlinear differential equations with Robin boundary, and the two dimensional Bratu-like equation [12]. In all these examples, the Newton–Raphson algorithm is employed to solve the corresponding nonlinear algebraic systems (36).

**Problem 1** Bratu equation

The Bratu equation is a typical example of nonlinear eigenvalue problems coming from various physical problems, such as the fuel ignition of thermal combustion and thermal reaction process in a rigid material [6], which in one-dimensional planar coordinates is of the form

$$\frac{d^2 u}{dx^2} + \lambda e^{u(x)} = 0, \quad u(0) = u(1) = 0. \tag{53}$$

For $0 < \lambda < \lambda_c \approx 3.51$, exact solution of Eq. (58) is $u(x) = -2\ln[\cosh(x\theta/2 - \theta/4)/\cosh(\theta/4)]$ with $\theta = \sqrt{2\lambda}\cosh(\theta/4)$ [6]. And for $\lambda = -1$, the exact solution becomes $u(x) = 2\ln[k\sec(kx/2 - k/4)] - \ln 2$.

Define $u_2 = d^2u/dx^2$, then Eq. (58) can be changed to

$$\begin{cases} u_2(x) + \lambda e^{u(x)} = 0 \\ u(x) = \int_0^x \int_0^\xi u_2(\xi_1)d\xi_1 d\xi - x\int_0^1 \int_0^\xi u_2(\xi_1)d\xi_1 d\xi \end{cases} \quad (54)$$

in which we have used $du/dx|_{x=0} = -\int_0^1 \int_0^\xi u_2(\xi_1)d\xi_1 d\xi$ determined by the boundary condition $u(1) = 0$. For a function $f(\cdot)$ and a vector $\mathbf{x} = \{x_k\}$, by defining $f.(\mathbf{x}) \triangleq \{f(x_k)\}$, then Eq. (59) can be discretized into algebraic equations by the proposed WICM as

$$\mathbf{u}_2 + \lambda e.^{\mathbf{A}_2 \mathbf{u}_2} \approx 0 \quad (55)$$

where $\mathbf{u}_2 = \{u_2(x_l), x_l = l/2^j\}$, $\mathbf{A}_2 = \{a_{2,k,l} = \Phi_{j,k}^{l_2}(x_l) - x_l \Phi_{j,k}^{l_2}(1)\}$, $k,l = 0,1,\cdots,2^j$. Once the $\mathbf{u}_2$ is obtained from Eq. (60), then nodal values of the unknown function can be obtained by $\mathbf{u} = \mathbf{A}_2 \mathbf{u}_2$.

Accuracy of the present WICM associated with the Coiflet with $N=6$, $M_1=7$ is compared with that of the WGM [6] and the SMBP [10]. Fig. 2 shows the absolute error at $x=1/2$ as a function of the number of grid points $n$, which is $2^j+1$ for the WICM, and $2^j-1$ for the WGM [6]. We can see from Fig. 2 that the convergence rate is about 7.0 for the proposed WICM and 4.0 for the WGM [6]. Fig. 3 shows the absolute value of the relative error at $x=1/2$ as a function of the parameter $\lambda$ under resolution level $j=6$, which is obtained by using both of the WICM and WGM [6]. It can be seen from Fig. 3 that the accuracy of WICM is much better than that of WGM [6] and is almost independent of the nonlinear intensity of the equation characterized by the parameter $\lambda$. Fig. 4(a) shows the distribution of the absolute error along $x$ by the WICM under resolution levels of $j=4, 5, 6$. Fig. 4(b) plots the maximal absolute error of the numerical results for $\lambda = -1$ as a function of the number of collocation points by the WICM and SMBP [10]. From Fig. 4(b), we can see that the error by WICM decays very fast as the resolution level, $j$, gets large. And from Fig. 4(b), the SMBP is found divergent as the number of collocation points or the degree of Bernstein polynomials increases, although the precision can be very high at a small polynomial degree. Such a drawback of the SMBP has also been recognized in other studies [29, 30].

**Problem 2** Nonlinear fourth-order differential equation

Consider the following nonlinear fourth-order differential equation

$$\begin{cases} u^{(4)}(x) - e^x u''(x) + u(x) + \sin(u(x)) = f(x), \quad 0 < x < 1 \\ f(x) = 1 + \sin(1 + \sinh(x)) - (e^x - 2)\sinh(x) \\ u(0) = 1, \ u'(0) = 1, \ u(1) = 1 + \sinh(1), \ u'(1) = \cosh(1) \end{cases} \quad (56)$$

Exact solution of Eq. (61) can be given by $u(x) = 1 + \sinh(x)$ as shown in [8]. By defining $u_i = d^i u / dx^i$ and applying boundary conditions, Eq. (61) can be equivalently changed to

$$\begin{cases} u_2 - e^x u_1 + u + \sin(u) = f(x) \\ u = u_4^{I_4}(x) + (2x^3 - 3x^2)u_4^{I_4}(1) + (x^2 - x^3)u_4^{I_3}(1) + ax^3/6 + bx^2/2 + xu'(0) + u(0) \\ u_2 = u_4^{I_2}(x) + (12x - 6)u_4^{I_4}(1) + (2 - 6x)u_4^{I_3}(1) + ax + b \end{cases}, (57)$$

where $a = -12u(1) + 6u'(1) + 12u(0) + 6u'(0)$ and $b = 6u(1) - 2u'(1) - 6u(0) - 4u'(0)$.

The final discretized algebraic system can be obtained as

$$\mathbf{u}_4 - \mathrm{diag}(\mathbf{B}_2 \mathbf{u}_2 + \mathbf{b}_2)e^{\mathbf{x}} + \mathbf{B}_4 \mathbf{u}_4 + \mathbf{b}_0 + \sin.(\mathbf{B}_4 \mathbf{u}_4 + \mathbf{b}_0) \approx f.(\mathbf{x}), \quad (58)$$

where matrices and vectors are

$$\mathbf{B}_2 = \{a_{2,k,l} = \Phi_{j,l}^{I_2}(x_k) + (12x_k - 6)\Phi_{j,l}^{I_4}(1) + (2 - 6x_k)\Phi_{j,l}^{I_3}(1)\},$$

$$\mathbf{B}_4 = \{a_{4,k,l} = \Phi_{j,l}^{I_4}(x_k) + (2x_k^3 - 3x_k^2)\Phi_{j,l}^{I_4}(1) + (x_k^2 - x_k^3)\Phi_{j,l}^{I_3}(1)\}, \mathbf{u}_i = \{u_i(x_l)\},$$

$$\mathbf{b}_0 = \{b_{0,k} = ax_k^3/6 + bx_k^2/2 + x_k u'(0) + u(0)\}, \ \mathbf{b}_2 = \{b_{2,k} = ax_k + b\}, \ x_l = l/2^j,$$

$$x_k = k/2^j \ k, l = 0, 1, \cdots, 2^j.$$

Numerical solutions are obtained by using the RKHSM [8], the new RKHSM (NRKSM) [9], the SMBP [10] and the proposed WICM, respectively. Absolute error of the numerical solutions obtained using these methods are listed in Table 4. We can see that the WICM with 17 collocation points has been much more accurate than those given by the RKHSM [8] and NRKSM [9] with 101, and the SMBP [10] with 17 grid points. Further, Fig. 5 illustrates the convergence of the SMBP and the WICM. The Bernstein series solution in terms of SMBP is again found not convergent with the degree of Bernstein polynomials.

**Problem 3** Fourth-order five-point BVP

We consider the Fourth-order five-point BVP as follows

$$\begin{cases} u^{(4)}(x) = \dfrac{u^{(2)}(x)}{1+u(x)} + q(x), \ 0 < x < 1 \\ q(x) = \theta^4 \sin(\theta x + \pi/6) + \theta^2 \sin(\theta x + \pi/6)/(1+\sin(\theta x + \pi/6)) \end{cases} \quad (59)$$

with boundary conditions

$$u(0) = \frac{1}{2}u(\frac{1}{2}), \quad u(1) = \frac{1}{4}u(\frac{1}{2}) + \frac{\sqrt{3}}{6}u(\frac{3}{4}), \quad u^{(2)}(0) = \frac{\sqrt{3}}{3}u^{(2)}(\frac{1}{4}),$$

$$u^{(2)}(1) = \frac{\sqrt{3}}{12}u^{(2)}(\frac{1}{4}) + \frac{\sqrt{3}}{4}u^{(2)}(\frac{3}{4}).$$

It can be easily verified that the exact solution to Eq. (64) is $u(x) = \sin(\theta x + \pi/6)$. Defining $u_i = d^i u / dx^i$ and applying the boundary conditions, Eq. (64) can be changed to

$$\begin{cases} u_4(x) = \dfrac{u_2(x)}{1+u(x)} + q(x) \\ u = u_4^{\int_4}(x) + \sum_{i=1}^{3}(\dfrac{x^3}{6}c_{4(i+1)} + \dfrac{x^2}{2}c_{3(i+1)} + xc_{2(i+1)} + c_{1(i+1)})u_4^{\int_2}(x_i) + \\ \sum_{j=1}^{3}(xc_{2j} + c_{1j})u_4^{\int_4}(x_j) \\ u_2 = u_4^{\int_2}(x) + \sum_{i=1}^{3}(xc_{4(i+3)} + c_{4(i+3)})u_4^{\int_2}(x_i) \end{cases} \quad (60)$$

By using the proposed WICM, we can obtain

$$(\mathbf{u}_4 - \mathbf{q}) + \text{diag}(\mathbf{B}_4\mathbf{u}_4)(\mathbf{u}_4 - \mathbf{q}) - \mathbf{B}_2\mathbf{u}_4 \approx 0 \quad (61)$$

where matrices and vectors are

$$\mathbf{u}_4 = \{u_4(x_l)\}, \ \mathbf{B}_2 = \{a_{2,k,l} = \Phi_{j,l}^{\int_2}(x_k) + \sum_{i=1}^{3}(x_k b_{3i} + b_{2i})\Phi_{j,l}^{\int_2}(x_i)\}, \text{ and}$$

$$\mathbf{B}_4 = \{a_{4,k,l} = \Phi_{j,l}^{\int_4}(x_k) + \sum_{i=1}^{3}(\dfrac{x_k^3}{6}b_{3i} + \dfrac{x_k^2}{2}b_{2i} + x_k b_{1i} + b_{0i})\Phi_{j,l}^{\int_2}(x_i) + \sum_{j=1}^{3}(x_k a_{1j} + a_{0j})\Phi_{j,l}^{\int_4}(x_j)\},$$

$$x_l = l/2^j, \ x_k = k/2^j \ k,l = 0,1,\cdots,2^j.$$

Numerical solutions are obtained from the fourth order CFDM [7] and the present WICM under number of grid points, $m=2^j+1$, where $j=4, 5, 6$ is the resolution levels. The maximal absolute error as a function of $m$ is plotted in Fig. 6(a), which shows that the convergence rate of the WICM reaches about 7.0. The absolute error distribution of the wavelet solutions under different resolution levels is given in Fig. 6(b), which shows that the absolute error of the WICM is about $O(10^{-10})$ when resolution level, $j = 4$, corresponding to $m=17$. By contrast, the solution by the CFDM with grid points $m = 128$ has a similar maximal absolute error $5.97\text{e}^{-10}$ [10]. Therefore, the proposed WICM can reach the same accuracy under much smaller computation cost than the CFDM.

**Problem 4** Large deformation bending of circular plates

The geometrically nonlinear equation of circular plates has been well studied both theoretically and numerically [32]. Choosing this problem as a benchmark is convenient because convergent analytical solutions for this problem have been found [32]. By define $\varphi$ and $S$ as the dimensionless deflection and internal axial force, such a problem can be described as [32]

$$\begin{aligned} y^2 \frac{d^2\varphi}{dy^2} &= \varphi(y)S(y) + y^2 Q \\ y^2 \frac{d^2 S}{dy^2} &= -\frac{1}{2}\varphi^2(y) \end{aligned}, \quad 0 \le y \le 1 \tag{62}$$

with Dirichlet and Robin boundary conditions

$$\varphi(0) = 0, \quad S(0) = 0; \quad \varphi(1) = \frac{\lambda}{\lambda - 1}\frac{d\varphi}{dy}\bigg|_{y=1}, \quad S(1) = \frac{\mu}{\mu - 1}\frac{dS}{dy}\bigg|_{y=1} \tag{63}$$

where $Q$ is the normalized uniform load applied on the circular plate, and $\lambda$ and $\mu$ are parameters determined by the supported condition along the boundary. For example, $\lambda$ and $\mu$ for the clamped plate are 0 and $2/(1-v)$, respectively, where $v = 0.3$ is the Poisson ratio. Further, define $\varphi_2 = d^2\varphi/dy^2$, $S_2 = d^2 S/dy^2$, the equivalent integral form of Eq. (67) can be derived as

$$\begin{cases} y^2 \varphi_2(y) = \varphi(y)S(y) + y^2 Q \\ y^2 S_2(y) = -\frac{1}{2}\varphi^2(y) \\ \varphi(y) = \varphi_2^{\int_2}(y) + y[(\lambda - 1)\varphi_2^{\int_2}(1) - \lambda\varphi_2^{\int}(1)] \\ S(y) = S_2^{\int_2}(y) + y[(\mu - 1)S_2^{\int_2}(1) - \mu S_2^{\int}(1)] \end{cases}. \tag{64}$$

Applying the WICM yields the final discretized algebraic system as follows:

$$\begin{aligned} \operatorname{diag}(\mathbf{y}^2)\boldsymbol{\varphi}_2 &\approx \operatorname{diag}(\mathbf{B}_2\mathbf{S}_2)(\mathbf{A}_2\boldsymbol{\varphi}_2) + Q\mathbf{y}^2 \\ \operatorname{diag}(\mathbf{y}^2)\mathbf{S}_2 &\approx -\frac{1}{2}\operatorname{diag}(\mathbf{A}_2\boldsymbol{\varphi}_2)(\mathbf{A}_2\boldsymbol{\varphi}_2) \end{aligned}, \tag{65}$$

where matrices and vectors are

$$\mathbf{A}_2 = \{a_{2,l,k} = \Phi_{j,k}^{\int_2}(y_l) + y_l(\lambda-1)\Phi_{j,k}^{\int_2}(1) - y_l\lambda\Phi_{j,k}^{\int}(1)\}, \quad \boldsymbol{\varphi}_2 = \{\varphi_2(y_l)\}, \mathbf{S}_2 = \{\varphi_2(y_l)\},$$

$$\mathbf{B}_2 = \{b_{2,l,k} = \Phi_{j,k}^{\int_2}(y_l) + y_l(\mu-1)\Phi_{j,k}^{\int_2}(1) - y_l\mu\Phi_{j,k}^{\int}(1)\}, \quad y_l = l/2^j, \text{ and } k,$$

$$l = 0, 1, \cdots, 2^j.$$

The entries of $a_{2,0,k}$ and $b_{2,0,k}$ are equal to zero, which causes singularity of the coefficient matrix. According to the L'Hôpital's rule and Eq. (67), we have

$$\varphi_2(0) = \left.\frac{d\varphi}{dy}\right|_{y=0} \left.\frac{dS}{dy}\right|_{y=0} + Q, \quad S_2(0) = -\frac{1}{2}\left(\left.\frac{d\varphi}{dy}\right|_{y=0}\right)^2. \tag{66}$$

Then the equations in Eq. (70) at $y_0 = 0$ can be replaced by

$$\varphi_2(0) \approx (\mathbf{a}_1\boldsymbol{\varphi}_2)(\mathbf{a}_2\mathbf{S}_2) + Q, \quad S_2(0) \approx -0.5(\mathbf{a}_1\boldsymbol{\varphi}_2)^2, \tag{67}$$

in which

$$\mathbf{a}_1 = \{a_{1,k} = (\lambda-1)\Phi_{j,k}^{\int_2}(1) - \lambda\Phi_{j,k}^{\int}(1)\} \text{ and } \mathbf{a}_2 = \{a_{2,k} = (\mu-1)\Phi_{j,k}^{\int_2}(1) - \mu\Phi_{j,k}^{\int}(1)\}.$$

Define the dimensionless central deflection as $W_m = -\int_0^1 \varphi(\xi)/\xi d\xi$. Wavelet solutions of the Eq. (67), in terms of the WGM [20] and the present WICM, are obtained under various resolution level, $j$, and normalized load $Q$. Fig. 7 shows the absolute value of the absolute error of $\varphi(1/2)$ and $S(1/2)$ as a function of the number of grid points $m=2^j+1$ under normalized load $Q=50$. We can see that the convergence rate of the WICM reaches about 7.0. Fig. 8 plots the distribution of absolute error on $\varphi$ and $S$ by the WICM with different resolution levels, which shows that the error level keeps similar at different locations, and rapidly decays as the resolution level $j$ increases. Fig. 9 plots the absolute value of the relative error on $W_m$ as a function of the normalized load $Q$ under resolution level $j=4$ for WICM and $j=6$ for WGM, respectively. It can be seen from Fig. (9) that the result by the WICM with only $2^4+1=17$ collocation points has been much more accuracy than that by the WGM [20] with $2^6-1=63$ grid points.

**Problem 5** Two dimensional Bratu-like equation

In order to verify if the proposed WICM can be applied to high dimensional NBVPs [12], we consider a typical example called the Bratu-like equation in two dimensions

$$\frac{\partial^2 u(x,y)}{\partial x^2}+\frac{\partial^2 u(x,y)}{\partial y^2}+\lambda e^{u(x,y)}=f(x,y),\ 0<x,y<1$$
$$u(x,y)|_{x=0,1}=0,\ u(x,y)|_{y=0,1}=0 \tag{68}$$

Exact solution of Eq. (73) is $u(x,y)=\sin(\pi x)\sin(\pi y)$ when $f(x,y)=\lambda e^{\sin(\pi x)\sin(\pi y)}-2\pi^2\sin(\pi x)\sin(\pi y)$. Defining $u_2=\partial^2 u/\partial x^2$, and $v_2=\partial^2 u/\partial y^2$, Eq. (73) can be changed to

$$\begin{cases} u_2(x,y)+v_2(x,y)+\lambda e^{u(x,y)}=f(x,y) \\ u(x,y)=u_2^{I_2}(x)-xu_2^{I_2}(1) \\ u(x,y)=v_2^{I_2}(y)-yv_2^{I_2}(1) \end{cases} \tag{69}$$

where we have considered $u_2^{I_2}(x)\triangleq\int_0^x\int_0^\xi u_2(\xi_1,y)d\xi_1 d\xi$ and $v_2^{I_2}(y)\triangleq\int_0^y\int_0^\eta v_2(x,\eta_1)d\eta_1 d\eta$. Applying Eq. (29) to approximate $u_2(x,y_{l'})$ and $v_2(x_{k'},y)$ for $l',k'=0,1,2,...,2^j$, we have

$$u_2(x,y_{l'})\approx\sum_{k=0}^{2^j}u_2(x_k,y_{l'})\Phi_{j,k}(x),\quad v_2(x_{k'},y)\approx\sum_{l=0}^{2^j}v_2(x_{k'},y_l)\Phi_{j,l}(y) \tag{75}$$

and

$$u(x,y_{l'})\approx\sum_{k=0}^{2^j}u_2(x_k,y_{l'})[\Phi_{j,k}^{I_2}(x)-x\Phi_{j,k}^{I_2}(1)],$$

$$v(x_{k'},y)\approx\sum_{l=0}^{2^j}v_2(x_{k'},y_l)[\Phi_{j,l}^{I_2}(y)-y\Phi_{j,l}^{I_2}(1)]. \tag{76}$$

Applying the WICM to Eq. (74) and considering Eqs. (75) and (76) leads to $\mathbf{u}\approx\mathbf{A}\mathbf{u}_2\approx\mathbf{B}\mathbf{v}_2$, and

$$\begin{cases} \mathbf{u}_2+\mathbf{v}_2+\lambda e^{\mathbf{A}\mathbf{u}_2}\approx\mathbf{f} \\ \mathbf{A}\mathbf{u}_2\approx\mathbf{B}\mathbf{v}_2 \end{cases} \tag{76}$$

where vectors and matrices are

$$\mathbf{u}_2=\{u_{2,K}=u_2(x_{k'},y_l)\},\ \mathbf{v}_2=\{v_{2,G}=v_2(x_k,y_{l'})\},\ \mathbf{f}=\{f_K=f(x_k,y_l)\},$$

$$\mathbf{A}=\{a_{G,K}=\Phi_{j,k}^{I_2}(x_{k'})-x_{k'}\Phi_{j,k}^{I_2}(1)\},\ \mathbf{B}=\{b_{G,K}=\Phi_{j,l}^{I_2}(y_{l'})-y_{l'}\Phi_{j,l}^{I_2}(1)\}$$

with subscripts $k,l,k',l'=0,1,\cdots,2^j$ and $K,G=1,...,(2^j+1)^2$. The vectors $\mathbf{u}_2$ and $\mathbf{v}_2$ can be solved from Eq. (76), then $\mathbf{u}$ can be obtained through $\mathbf{u}\approx\mathbf{A}\mathbf{u}_2$ or $\mathbf{u}\approx\mathbf{B}\mathbf{v}_2$ in terms of Eq. (76).

Fig. 10(a) shows the absolute error on $u(1/2, 1/2)$ as a function of $m=2^j$. We note that the number of grid points along one dimension for the WICM is $2^j+1$, and for the WGM [12] is $2^j-1$. We can see that the WICM has a convergence rate, 8.0, which is much higher than that of the WGM [12], 5.0. Absolute error on $u(1/2, 1/2)$ as a function of parameter $\lambda$ is ploted in Fig. 10(b), which, interestingly, shows that the absolute error of the WICM is almost unchanged for a wide range of $\lambda$, which characterizes the nonlinearity intensity of Eq. (73), implying that the proposed WICM is suitable for BVPs with strong nonlinearity.

## 5. Conclusion

In the present study, we have proposed the WICM for general BVPs. It has been demonstrated both analytically and numerically that this method can reach an accuracy of $O(2^{-jN})$, where $j$ is the resolution level and $N$-1 corresponds to the number of vanishing moments of the adopted Coiflet. As the Coiflet can be arbitrarily chosen, allowing $N$ to be any positive even numbers, therefore the proposed WICM can theoretically have arbitrarily high order of accuracy. Most interestingly, unlike most other existing numerical methods, accuracy of the WICM is independent of the order of the equation to be solved. And condition numbers of corresponding matrices after spatial discretization almost keep independent of the number of grid points associated with the resolution level $j$. Typical numerical examples have justified the accuracy, efficiency and robustness of the proposed method.

**Acknowledgement**

This research is supported by grants from the National Natural Science Foundation of China (11421062, 11502103).

**Appendix A**

We have defined the $n$-tuple integral of a scaling function $\phi(2^j x - l)$, $l \in \mathbf{Z}$ as

$$\phi_{j,l}^{\int_n}(x) = \int_0^x \int_0^{\xi_n} \cdots \int_0^{\xi_2} \phi(2^j \xi_1 - l) d\xi_1 d\xi_2 \cdots d\xi_n$$
$$= \frac{1}{2^{jn}} \phi^{\int_n}(2^j x - l) - \frac{1}{2^{jn}} \sum_{m=1}^{n} \frac{(2^j x)^{n-m}}{(n-m)!} \phi^{\int_m}(-l) , \quad \text{(A1)}$$

In order to obtain the value of $\phi^{\int_n}(x)$, by integrating Eq. (4) we have

$$\phi^{\int_n}(x) = 2^{-n}\sum_{k=0}^{3N-1} p_k \phi^{\int_n}(2x-k). \tag{A2}$$

Noting that $\phi^{\int_1}(3N-1)=1$, and $\phi^{\int_1}(x\leq 0)=0$, we can easily derive the following properties of $\phi^{\int_n}(x)$ [25]

$$\phi^{\int_n}(x) = 0, \quad x \leq 0, \tag{A3}$$

$$\phi^{\int_n}(x) = \sum_{k=0}^{n-1} \frac{(x-3N+1)^k}{k!}\phi^{\int_{n-k}}(3N-1), \quad x \geq 3N-1, \tag{A4}$$

$$\phi^{\int_n}(3N-1) = \frac{1}{2^n-2}\sum_{j=1}^{n-1}\sum_{k=0}^{3N-1} p_k \frac{(3N-1-k)^j}{j!}\phi^{\int_{n-j}}(3N-1). \tag{A5}$$

Based on Eqs. (A2-A5), the values of $\phi^{\int_n}(x)$ for $x=1,2,\cdots,3N-2$ can be determined from the following linear equations [25]

$$(\mathbf{I}-2^{-n}\mathbf{M})\mathbf{\Phi}^{\int_n} = \mathbf{c} \tag{A6}$$

where $\mathbf{I}$ is the unit matrix, and

$$\mathbf{M} = \begin{bmatrix} p_1 & p_0 & \cdots & 0 & 0 \\ p_3 & p_2 & \cdots & 0 & 0 \\ \vdots & \vdots & \vdots & \vdots & \vdots \\ 0 & 0 & \cdots & p_{3N-3} & p_{3N-4} \\ 0 & 0 & \cdots & p_{3N-1} & p_{3N-2} \end{bmatrix}, \tag{A7}$$

$$\mathbf{\Phi}^{\int_n} = \{\phi^{\int_n}(1),\phi^{\int_n}(2),\cdots,\phi^{\int_n}(3N-2)\}^T \text{ and } \mathbf{c} = \{c_i = \sum_{\substack{2i-k\geq 3N-1 \\ k=0,1,2,\cdots,3N-1}} p_k \phi^{\int_n}(2i-k)\}^T. \tag{A8}$$

## Appendix B

From Eq. (30), we have

$$\Phi_{j,k}^{\int_i}(\frac{l}{2^j}) = \begin{cases} \phi_{j,k-M_1}^{\int_i}(\frac{l}{2^j}) + \sum_{m=-\alpha_2}^{0} T_{L,k}(\frac{m}{2^j})\phi_{j,m-M_1}^{\int_i}(\frac{l}{2^j}), & 0 \leq k \leq \alpha_1 \\ \phi_{j,k-M_1}^{\int_i}(\frac{l}{2^j}), & \alpha_1+1 \leq k \leq 2^j-\alpha_2-1 \\ \phi_{j,k-M_1}^{\int_i}(\frac{l}{2^j}) + \sum_{m=1}^{\alpha_1} T_{R,2^j-k}(1+\frac{m}{2^j})\phi_{j,2^j+m-M_1}^{\int_i}(\frac{l}{2^j}), & 2^j-\alpha_2 \leq k \leq 2^j \end{cases}, \tag{B1}$$

where

$$T_{L,k}(\frac{m}{2^j}) = \sum_{i=0}^{N-1} 2^{ij}\frac{\zeta_{a,i,k}}{i!}\frac{m^i}{2^{ij}} = \sum_{i=0}^{N-1} m^i \frac{\zeta_{a,i,k}}{i!}, \quad m=-\alpha_2,\ldots,0, \tag{B2}$$

$$T_{R,2^j-k}(1+\frac{m}{2^j}) = \sum_{i=0}^{N-1} 2^{ij}\frac{\zeta_{b,i,2^j-k}}{i!}\frac{m^i}{2^{ij}} = \sum_{i=0}^{N-1}\frac{\zeta_{b,i,2^j-k}}{i!}m^i, \quad m=1,\ldots,\alpha_1, \qquad (B3)$$

Eqs. (B2) and (B3) show that $T_{L,k}(\frac{m}{2^j})$ and $T_{R,2^j-k}(1+\frac{m}{2^j})$ are actually independent of the resolution level $j$. And from Eq. (A1), we have

$$\phi_{j,k-M_1}^{\int_i}(\frac{l}{2^j}) = \frac{1}{2^j}\Lambda_{i,l,k} \qquad (B4)$$

where

$$\Lambda_{i,l,k} = \frac{1}{2^{(i-1)j}}\phi^{\int_i}(l-k+M_1) - \sum_{m=1}^{i}\frac{l^{i-m}}{(i-m)!}\frac{1}{2^{(i-1)j}}\phi^{\int_m}(M_1-k). \qquad (B5)$$

Considering Eqs. (A3)-(A5) and noting $l, k=0, 1, \ldots, 2^j$, we can easily verify that $\Lambda_{i,l,k} \to$ constant as $j$ approaching infinity, so as to $\phi_{j,k-M_1}^{\int_i}(\frac{l}{2^j}) = O(\frac{1}{2^j})$. Submitting Eq. (B4) into Eq. (B1) gives the expression of $\Phi_{j,k}^{\int_i}(\frac{l}{2^j})$ on the values of $\phi^{\int_i}(x)$ at integer points. If we further define

$$\mathrm{K}_{i,l,k} = \begin{cases} \Lambda_{i,l,k} + \sum_{m=-\alpha_2}^{0} T_{L,k}(\frac{m}{2^j})\Lambda_{i,l,m}, & 0 \leq k \leq \alpha_1 \\ \Lambda_{i,l,k}, & \alpha_1+1 \leq k \leq 2^j-\alpha_2-1, \\ \Lambda_{i,l,k} + \sum_{m=1}^{\alpha_1} T_{R,2^j-k}(1+\frac{m}{2^j})\Lambda_{i,l,2^j+m}, & 2^j-\alpha_2 \leq k \leq 2^j \end{cases} \qquad (B6)$$

then Eq. (B1) can be rewritten as

$$\Phi_{j,k}^{\int_i}(\frac{l}{2^j}) = \frac{1}{2^j}\mathrm{K}_{i,l,k} \qquad (B7)$$

with $\mathrm{K}_{i,l,k}$ is on the order of $O(1)$.

Table 1. Filter coefficients $p_k$ for the Coiflet with $N=6$ and $M_1=7$.

| k | $p_k$ | k | $p_k$ | K | $p_k$ |
|---|---|---|---|---|---|
| 0 | -2.3926386572801E-03 | 6 | 6.4599454329399E-01 | 12 | 1.2388695657060E-02 |
| 1 | -4.9326018541804E-03 | 7 | 1.1162662132580E+00 | 13 | -1.5831780392559E-02 |
| 2 | 2.7140399711400E-02 | 8 | 5.3818905570800E-01 | 14 | -2.7171786005400E-03 |
| 3 | 3.0647555946200E-02 | 9 | -9.9615433862400E-02 | 15 | 2.8869486640200E-03 |
| 4 | -1.3931023707080E-01 | 10 | -7.9923139434800E-02 | 16 | 6.3049939470800E-04 |
| 5 | -8.0606530717800E-02 | 11 | 5.1491462932400E-02 | 17 | -3.0583397359600E-04 |

Table 2. The values of $\phi^{j_n}(k)$ for $N=6$, $M_1=7$, $n=1, 2, 3$ and 4.

| k | $\phi^j(k)$ | $\phi^{j_2}(k)$ | $\phi^{j_3}(k)$ | $\phi^{j_4}(k)$ |
|---|---|---|---|---|
| 1 | -4.527262019294291Ee-09 | -2.640542520280718E-10 | -1.177917982279285E-11 | -4.830484119726849E-13 |
| 2 | 3.793659018256384E-06 | 4.419887973421076E-07 | 3.940901828183616E-08 | 3.231226434268301E-09 |
| 3 | 2.929388903777156E-06 | -3.394779202257466E-06 | -1.601106900878482E-06 | -1.390327115561728E-07 |
| 4 | -3.134173582887661E-03 | -7.269055985025437E-04 | -1.280198775311740E-04 | -2.128451273233253E-05 |
| 5 | 2.264311126637817E-02 | 3.381871078932003E-03 | -1.008794458421603E-03 | -8.249450856216146E-04 |
| 6 | -8.486417521139951E-02 | -9.611317252824955E-03 | 5.958177925875930E-03 | 2.387016412381357E-03 |
| 7 | 5.234649707003830E-01 | 9.430655504920177E-02 | -3.787403971467576E-03 | -3.184751382798324E-03 |
| 8 | 1.051493179855499E+00 | 9.976688497344420E-01 | 4.980359246292250E-01 | 1.666535756413486E-01 |
| 9 | 9.875787622086671E-01 | 1.997725744751554E+00 | 2.001481573098073E+00 | 1.333433158829295E+00 |
| 10 | 1.003591003107596E+00 | 3.000516645305960E+00 | 4.499416305495880E+00 | 4.500187806094493E+00 |
| 11 | 9.991452076341547E-01 | 4.000084460530419E+00 | 8.000033233268191E+00 | 1.066664726952882E+01 |
| 12 | 1.000073554299164E+00 | 4.999990068904882E+00 | 1.250000065760564E+01 | 2.083333339909631E+01 |
| 13 | 1.000001685301915E+00 | 6.000000321738184E+00 | 1.799999990934045E+01 | 3.600000000204158E+01 |
| 14 | 1.000000156123390E+00 | 6.999999992224733E+00 | 2.449999999879408E+01 | 5.716666666713520E+01 |
| 15 | 9.999999997921172E-01 | 7.999999999985318E+00 | 3.200000000005910E+01 | 8.533333333347005E+01 |
| 16 | 1.000000000000944E+00 | 9.000000000012078E+00 | 4.050000000006736E+01 | 1.215000000001882E+02 |
| 17 | 1.000000000000000E+00 | 1.000000000000843E+01 | 5.000000000006850E+01 | 1.666666666668998E+02 |

Table 3. Condition number of matrices $\mathbf{A}_i$ and norm of their inverse.

| $J$ | $\left\|\mathbf{A}_1^{-1}\right\|_2$ | $\left\|\mathbf{A}_2^{-1}\right\|_2$ | $\left\|\mathbf{A}_3^{-1}\right\|_2$ | $K_2(\mathbf{A}_1)$ | $K_2(\mathbf{A}_2)$ | $K_2(\mathbf{A}_3)$ |
|---|---|---|---|---|---|---|
| 4 | 1.0261 | 2.1037 | 1.0135 | 1.4944 | 2.1130 | 1.0552 |
| 5 | 1.0135 | 2.0536 | 1.0127 | 1.4661 | 2.0404 | 1.0530 |
| 6 | 1.0069 | 2.0274 | 1.0123 | 1.4516 | 2.0028 | 1.0519 |
| 7 | 1.0035 | 2.0139 | 1.0121 | 1.4442 | 1.9838 | 1.0513 |

Table 4. Absolute error of the numerical solutions of Eq. (61).

| $x$ | RHKSM, $m=101$ [8] | NRKSM, $m=101$ [9] | SMBP, $m=17$ [10] | WICM, $m=17$ |
|---|---|---|---|---|
| 0.1 | 2.78e-8 | 7.5e-9 | 2.92411e-12 | 1.33227e-15 |
| 0.2 | 8.09e-8 | 1.5e-9 | 4.00524e-12 | 3.10862e-15 |
| 0.3 | 1.20e-7 | 2.1e-9 | 6.16573e-12 | 5.99520e-15 |
| 0.4 | 1.25e-7 | 2.3e-9 | 8.32978e-12 | 7.32747e-15 |
| 0.5 | 9.56e-8 | 2.4e-9 | 9.41269e-12 | 6.66134e-15 |
| 0.6 | 4.82e-8 | 2.4e-9 | 1.15699e-11 | 3.99680e-15 |
| 0.7 | 7.38e-9 | 2.4e-9 | 1.26406e-11 | 1.99840e-15 |
| 0.8 | 1.07e-8 | 2.2e-9 | 1.48117e-11 | 1.11022e-15 |
| 0.9 | 7.08e-9 | 1.4e-9 | 1.56586e-11 | 1.55431e-15 |

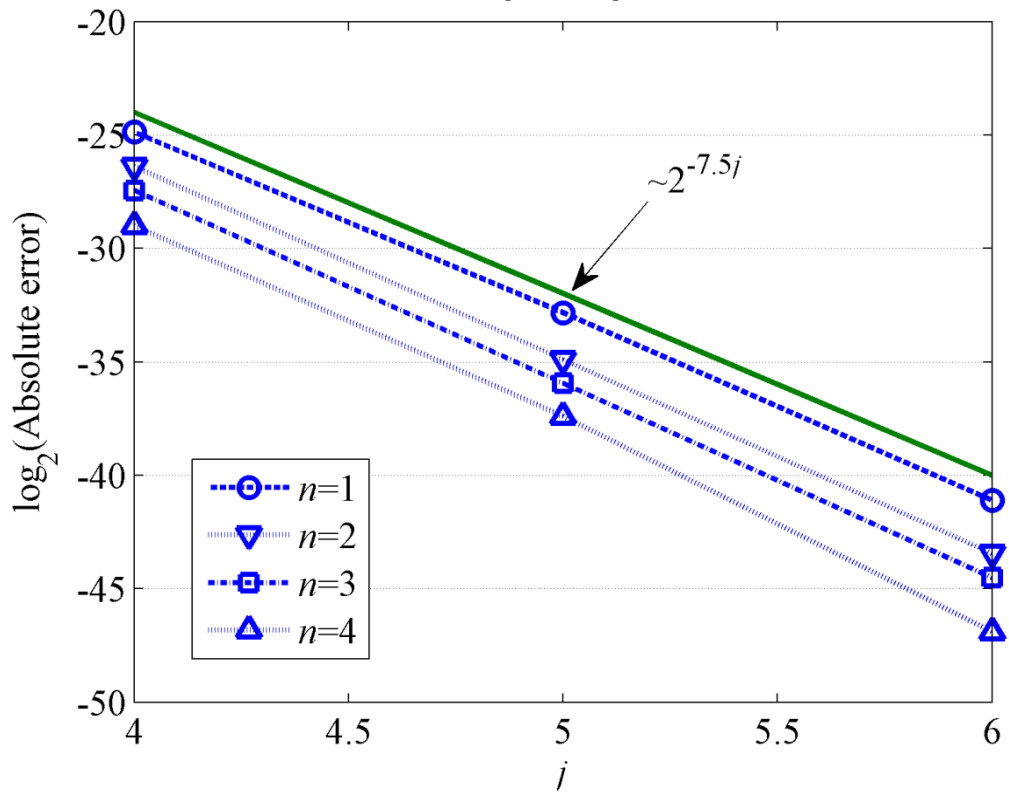

Fig. 1. Absolute error for the wavelet approximation of *n*-tuple integrals of $\sin(\pi x)$ as a function of resolution *j*, where *n* = 1, 2, 3, 4.

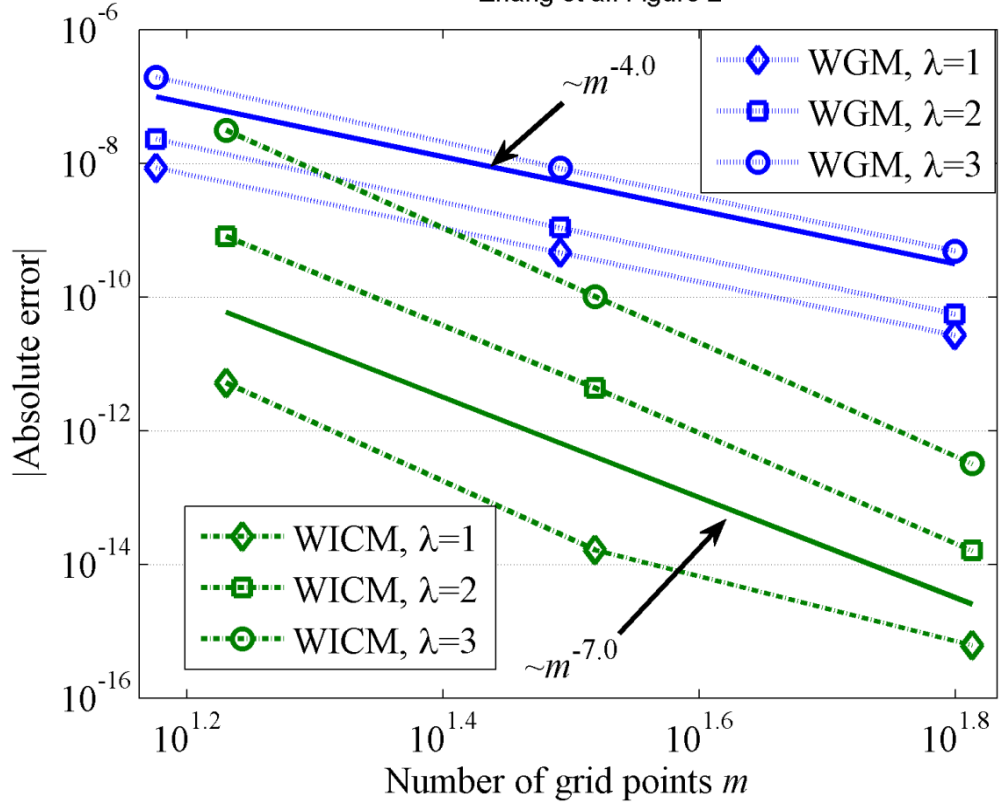

Fig. 2. Absolute error at $x=1/2$ as a function of the number of grid points $m=2^j+1$ for the WICM and $m=2^j-1$ for the WGM [6]

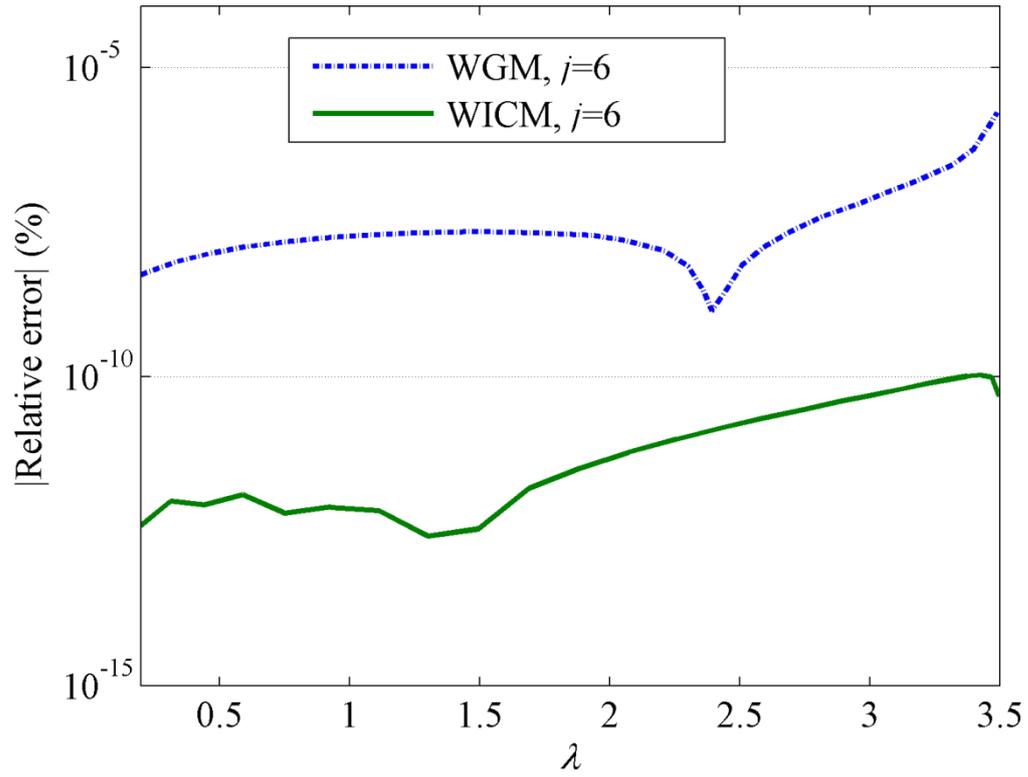

Fig. 3. Absolute values of relative errors of numerical results of $u(1/2)$ as a function of parameter $\lambda$ under resolution level $j$=6.

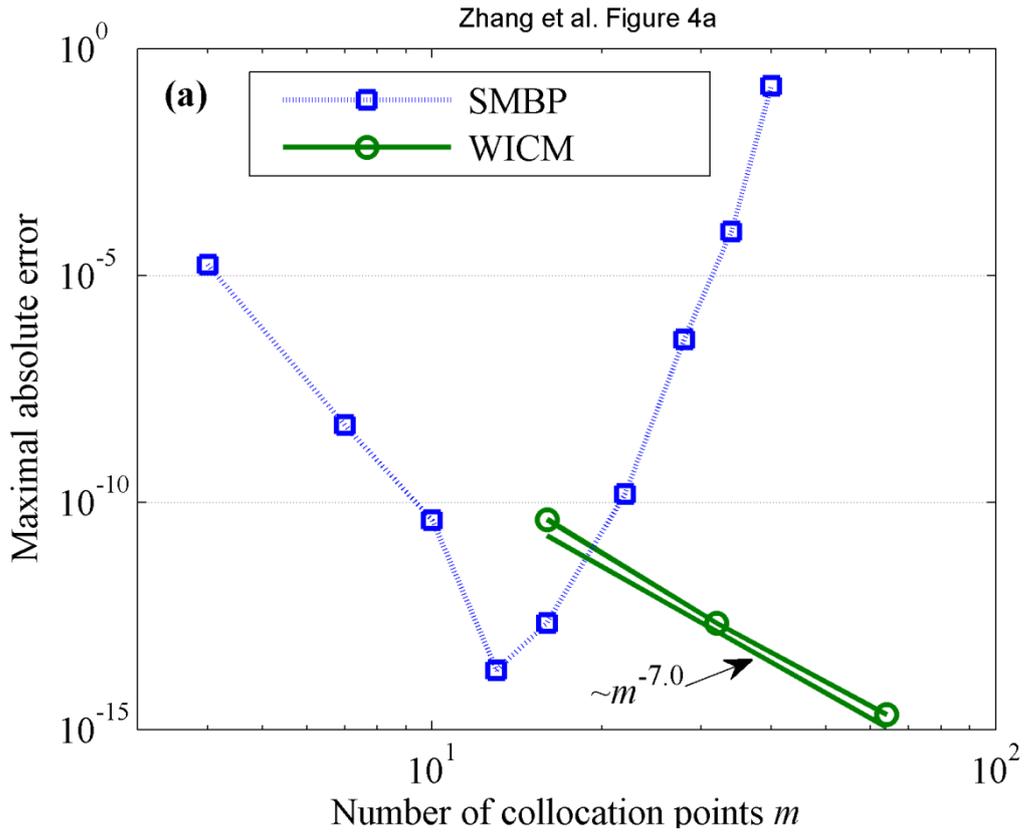

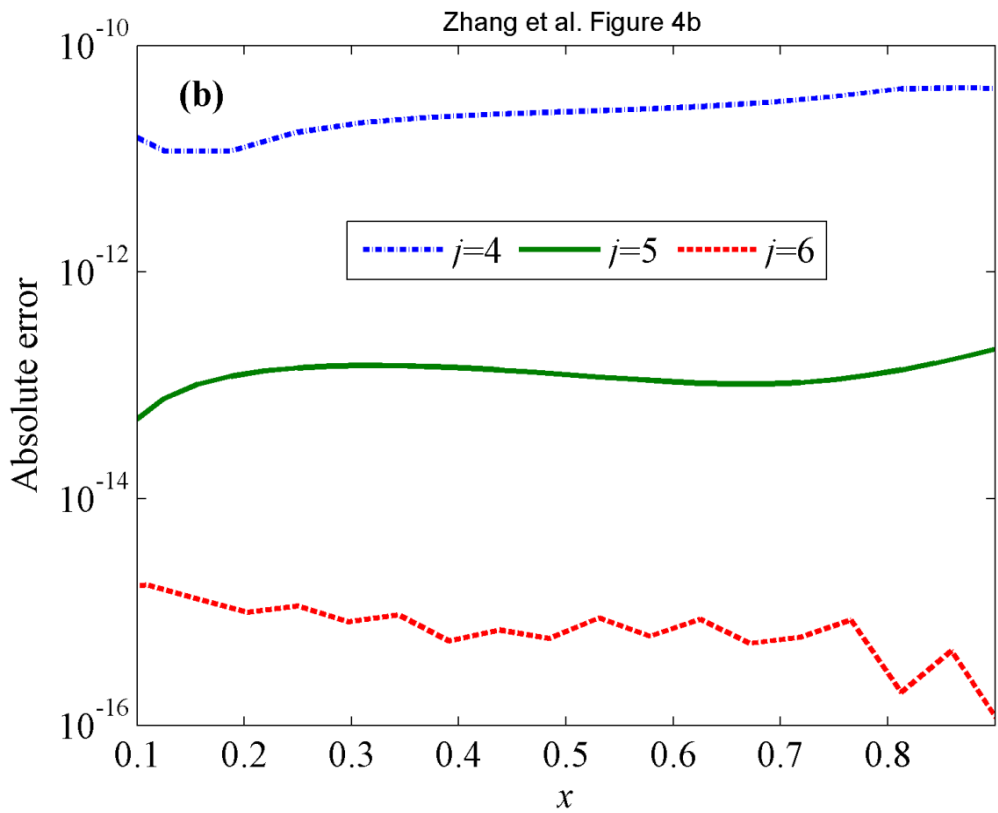

Fig. 4. Absolute error of numerical solutions of Eq. (58) with $\lambda=-1$: (a) Distribution of the absolute error by WICM for $j=4, 5, 6$; (b) Maximal absolute error as a function of number of collocation points, $m=2^j+1$.

\

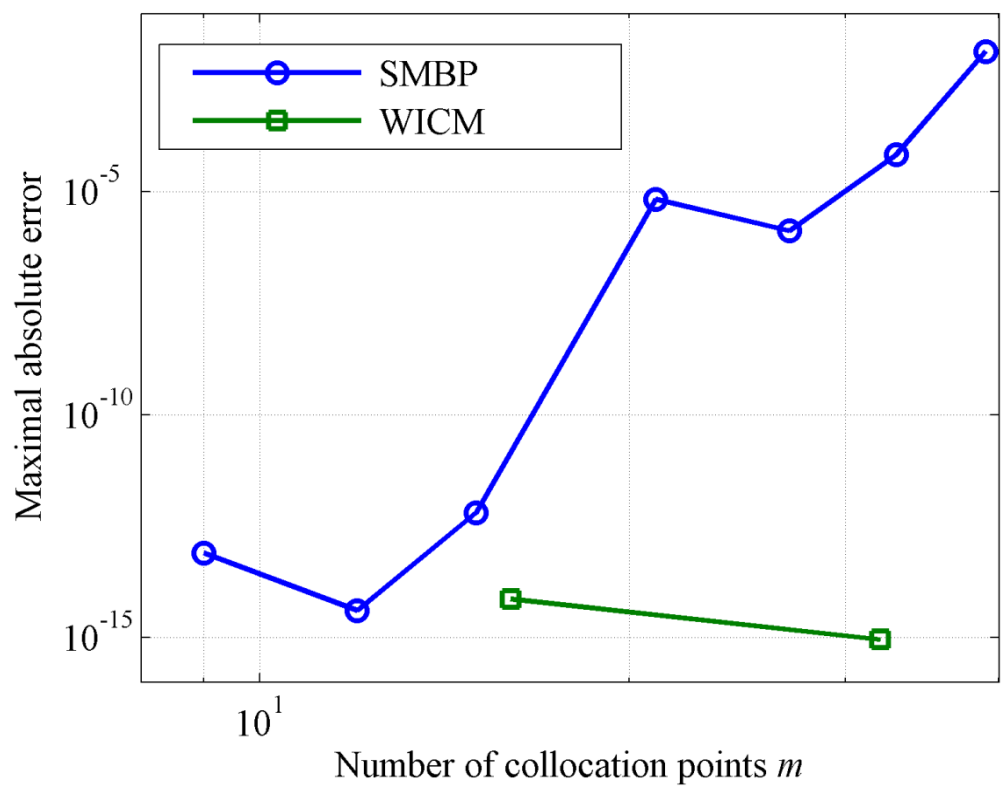

Fig. 5. Maximal absolute errors as a function of the number of collocation points $m=2^j+1$ with $j$ the resolution level.

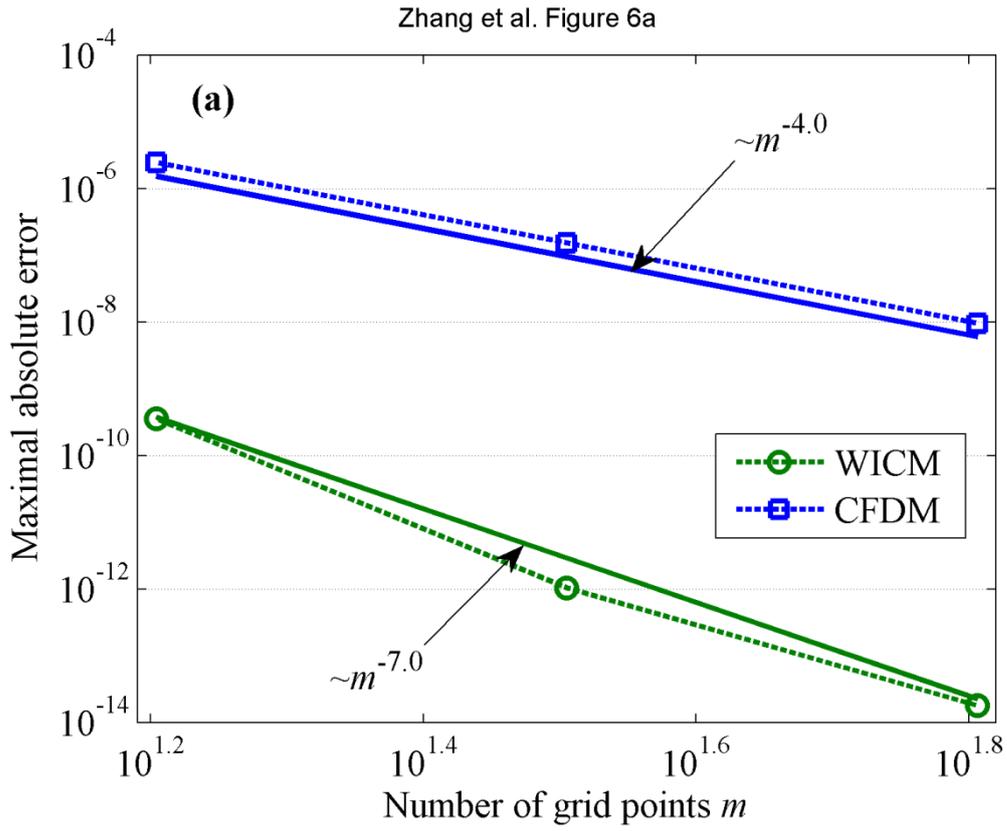

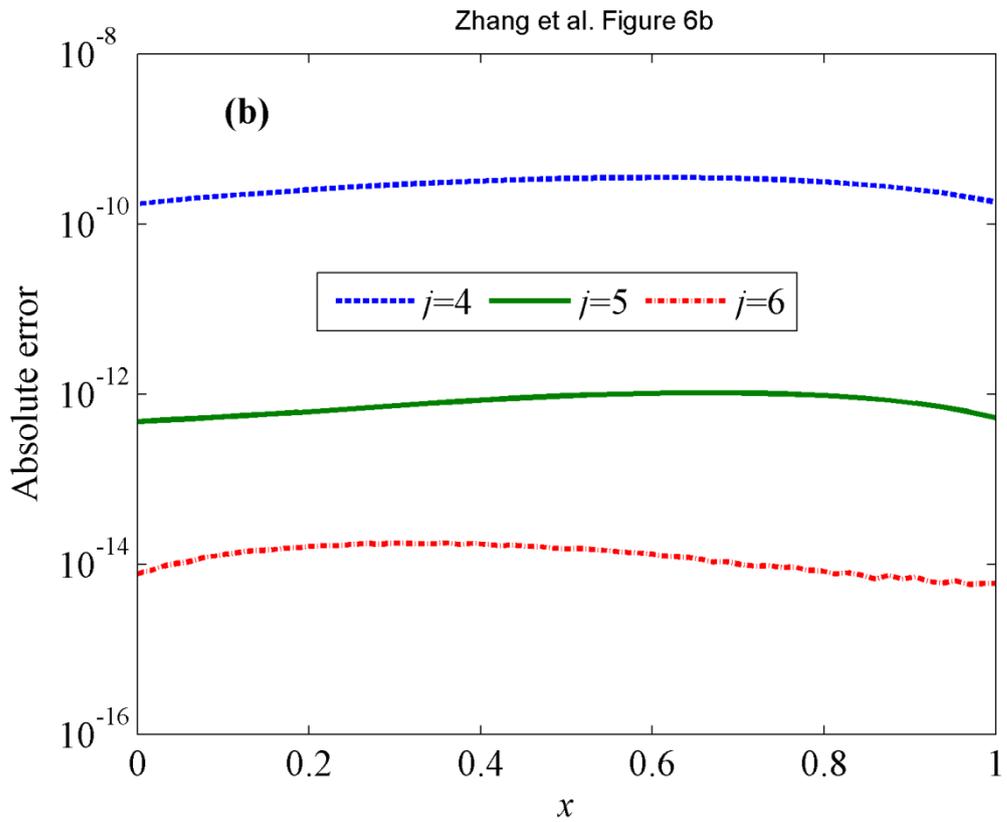

Fig. 6. Absolute error of numerical solutions of Eq. (64): (a) Dependence of the Maximal absolute error on the number of grid points $m=2^j+1$ with $j$ the resolution level; (b) Distribution of the absolute error along $x$ by the WICM.

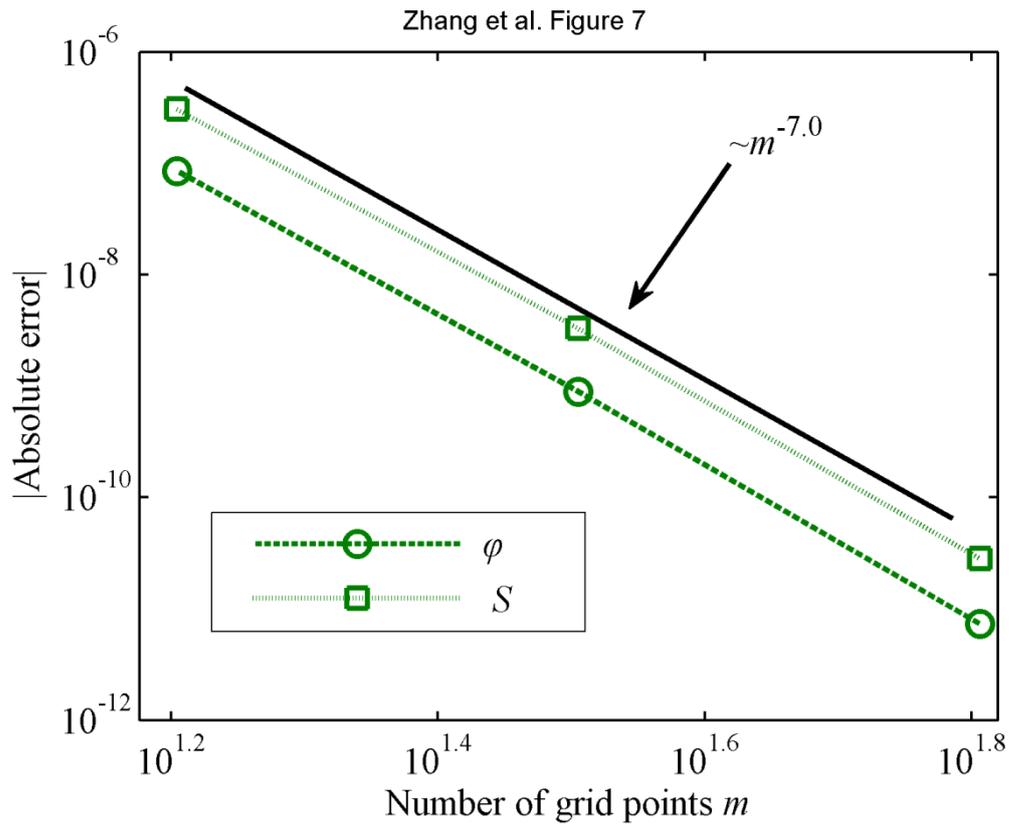

Fig.7. Absolute value of the absolute error at $x=1/2$ as a function of the number of grid points $m=2^j+1$ under normalized load $Q=50$, where $j$ is the resolution level, $\varphi$ and $S$ are dimensionless deflection and internal axial force.

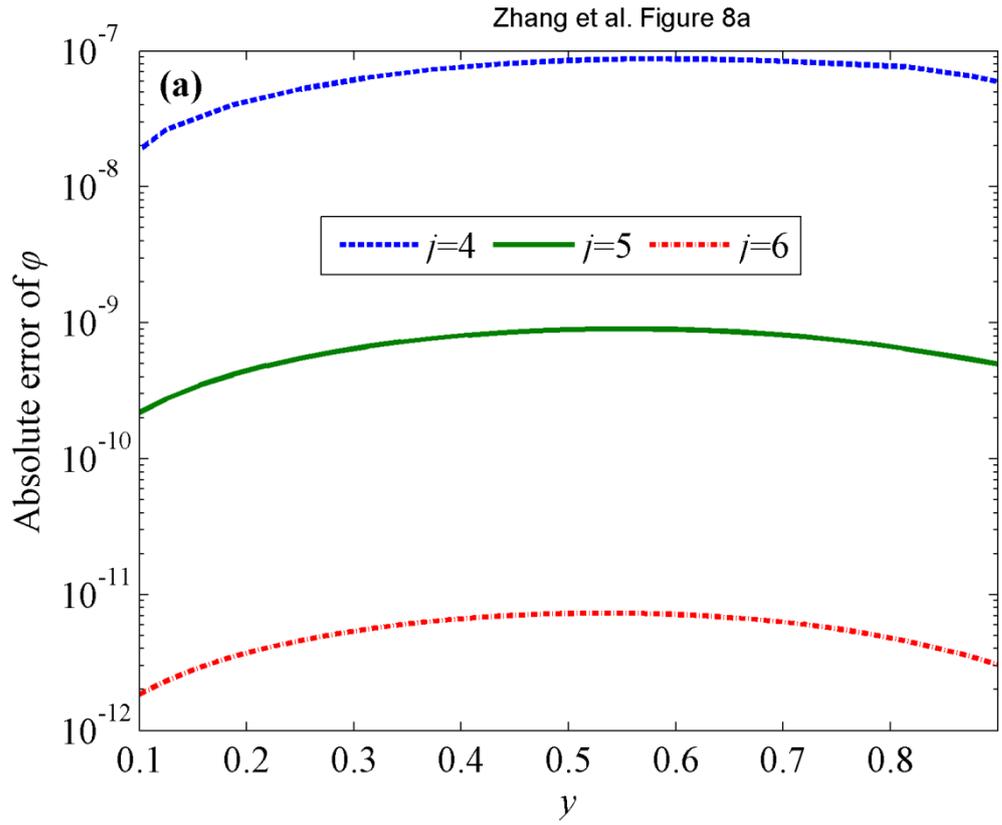

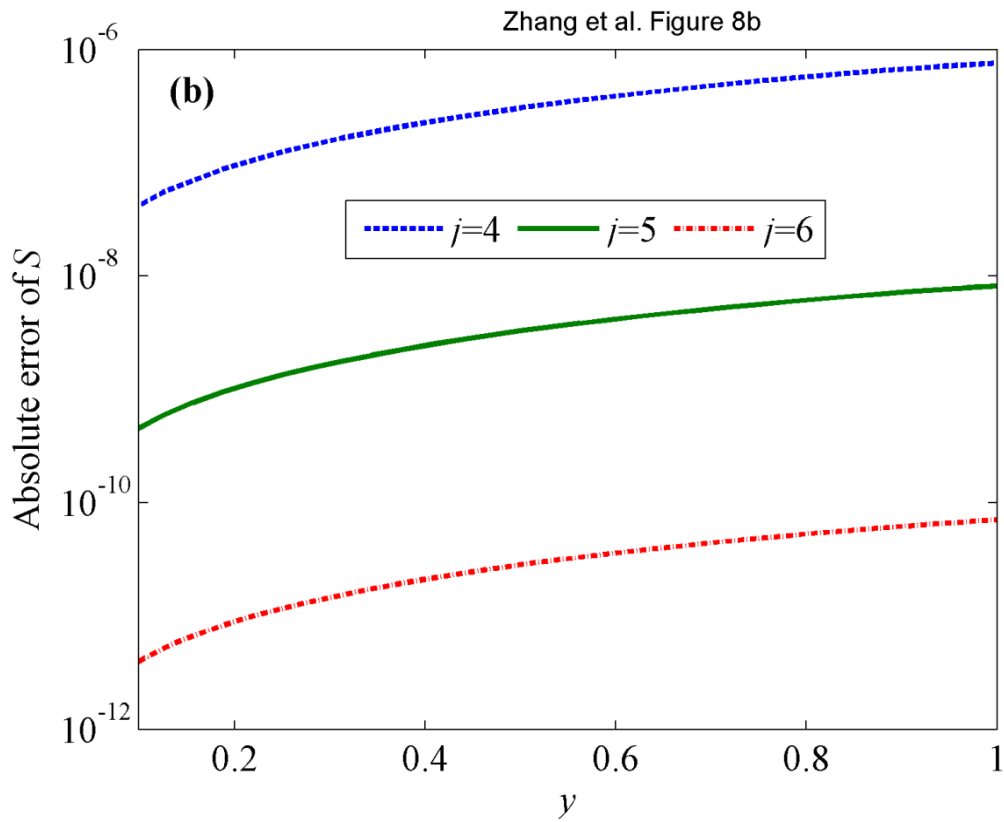

Fig. 8. Error distribution of the numerical solutions of Eq. (67) by the WICM, where $Q$=50, resolution level $j$=4, 5, 6: (a) Distribution of the absolute error of $\varphi$ along $y$; (b) Distribution of the absolute errors of $S$ along $y$.

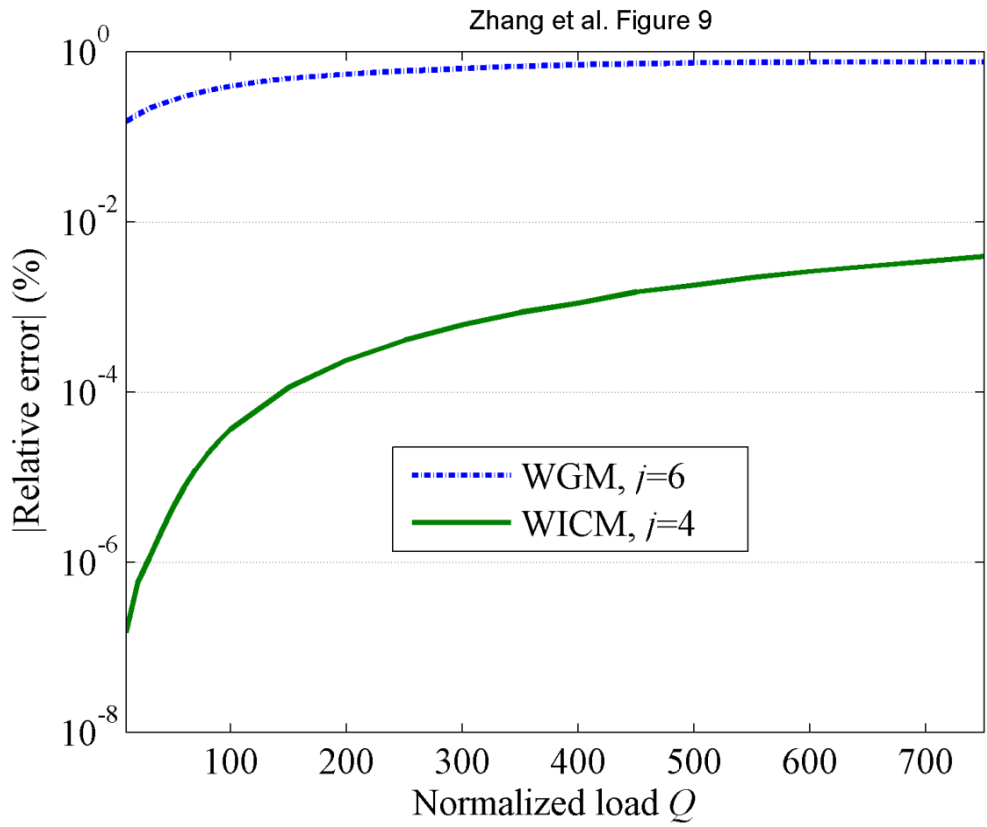

Fig. 9. Absolute value of the relative error of numerically obtained $W_m$ as a function of the normalized load $Q$ under resolution level $j=4$ for WGM [20] and $j=6$ for WICM, respectively.

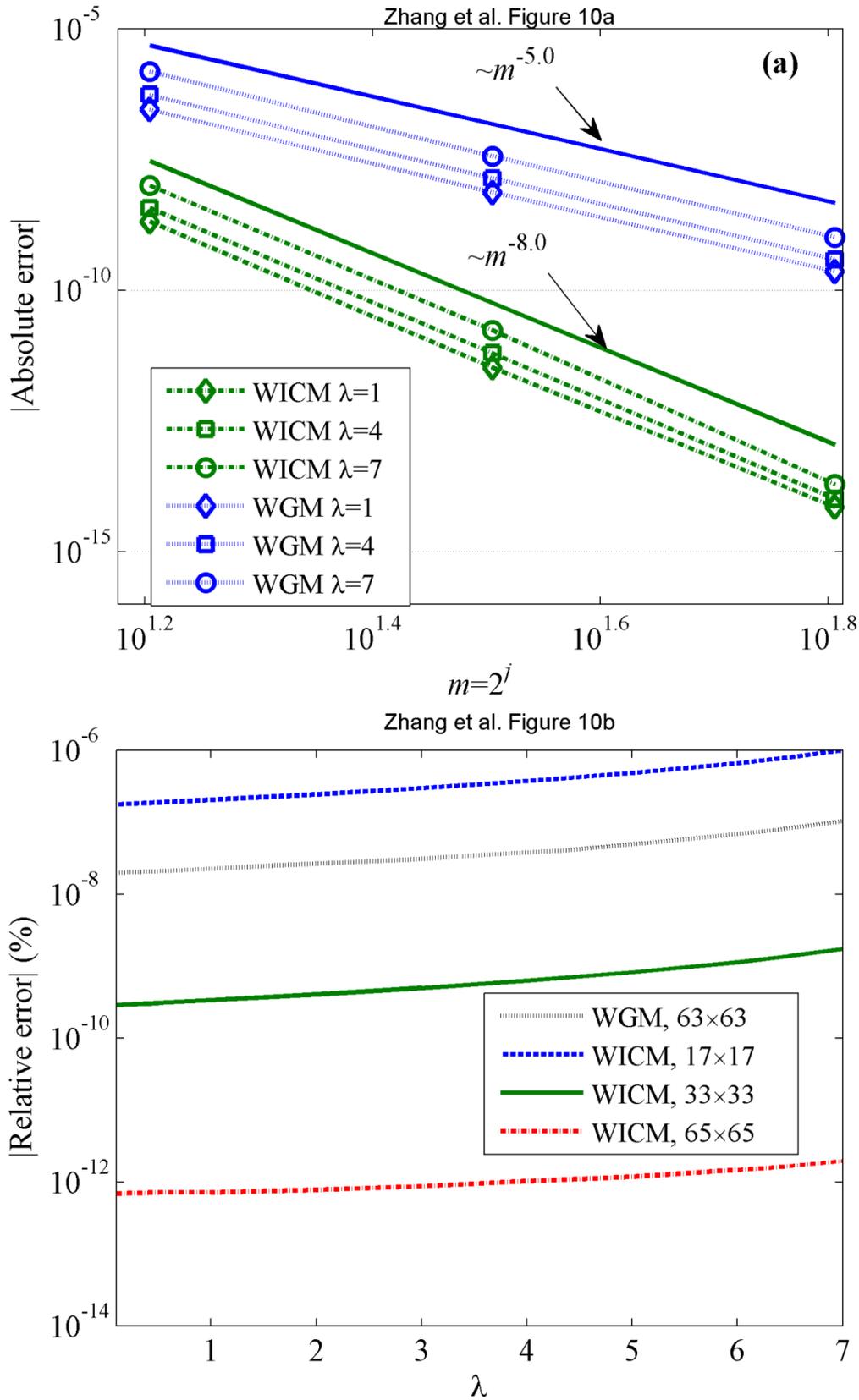

Fig. 10. Error of the numerical solution of Eq. (73) under different parameters, $\lambda$, and resolution levels, $j$, by the WGM [12] and WICM: (a) Absolute error as a function of the number of grid points $m=2^j$ along one dimension; (b) Relative error as a function

of parameter $\lambda$ under different grid number